\documentclass[11pt, twoside]{article}
\usepackage{latexsym}
\usepackage{amsmath}
\usepackage{amssymb}
\usepackage[all]{xy}
\usepackage{amsfonts}
\usepackage{verbatim}
\usepackage{amsthm}
\usepackage{mathrsfs}
\usepackage{epsfig}
\usepackage{xy}
\usepackage{array}
\usepackage{stmaryrd}
\usepackage{graphicx,color}
\usepackage{xcolor}
\usepackage[colorlinks=true,linkcolor=blue,citecolor=blue]{hyperref}
\usepackage{tikz}
\usetikzlibrary{arrows,calc}
\usepackage{etex}
\usepackage{mathdots}
\usepackage{float}
\usepackage{graphics}
\usepackage{pdflscape}

\usepackage{anysize}
\input xypic
\xyoption{all}
\usepackage{enumitem}
\usepackage[perpage,symbol]{footmisc}
\usepackage{geometry}
\usepackage{setspace}
\def\A{\mathcal{A}}

\def\C{\mathcal{C}}
\def\E{\mathbb{E}}

\def\dr{\ar@{->}[r]}

\numberwithin{equation}{section}

\newcommand{\add}{\mathsf{add}\hspace{.01in}}

\def\proj{\mathsf{proj}\hspace{.01in}}
\newcommand{\coker}{\operatorname{coker}}
\def\End{\mathsf{End}}
\def\Ext{\mbox{Ext}}
\def\Hom{\mbox{Hom}}

\usepackage{bm}
\begin{document}
\baselineskip=15pt

\def\Pcal{\mathcal{P}}
\def\Ical{\mathcal{I}}
\newcommand{\ind}{\operatorname{ind}}
\newcommand{\Supp}{\operatorname{Supp}}
\newcommand{\rad}{\operatorname{rad}}
\newcommand{\AR}{\operatorname{AR}}
\newcommand{\Rel}{\operatorname{Rel}}
\newcommand{\Proj}{\operatorname{Proj}}
\newcommand{\Inj}{\operatorname{Inj}}
\newcommand{\Db}{\mathrm D^{\mathrm b}}
\newcommand{\modcat}{\operatorname{mod}}
\newcommand{\length}{\operatorname{length}}

\title{\Large{\bf Auslander--Reiten relations and local finiteness in\\[2mm] extriangulated categories}}
\medskip
\author{Panyue Zhou and Bin Zhu}
\date{}

\maketitle

\newtheorem{theorem}{Theorem}[section]
\newtheorem{lemma}[theorem]{Lemma}
\newtheorem{corollary}[theorem]{Corollary}
\newtheorem{proposition}[theorem]{Proposition}
\newtheorem{conjecture}{Conjecture}
\theoremstyle{definition}
\newtheorem{definition}[theorem]{Definition}
\newtheorem{question}[theorem]{Question}
\newtheorem{remark}[theorem]{Remark}
\newtheorem{remark*}[]{Remark}
\newtheorem{example}[theorem]{Example}
\newtheorem{example*}[]{Example}
\newtheorem{condition}[theorem]{Condition}
\newtheorem{condition*}[]{Condition}
\newtheorem{construction}[theorem]{Construction}
\newtheorem{construction*}[]{Construction}
\newtheorem{assumption}[theorem]{Assumption}
\newtheorem{assumption*}[]{Assumption}

\baselineskip=17pt
\parindent=0.5cm
\vspace{-6mm}

\begin{abstract}
\baselineskip=16pt
Let $\C$ be a Hom-finite Krull--Schmidt extriangulated category over a commutative artin ring.
Write $\Rel(\C)$ for the subgroup of the split Grothendieck group generated by the
relations of all $\E$-triangles, and $\AR(\C)$ for the subgroup generated by
those coming from Auslander--Reiten $\E$-triangles.  We first prove that either one-sided finiteness condition on Hom supports implies
$\Rel(\C)=\AR(\C)$.
For incoming Hom supports, the proof uses a Jordan--H\"{o}lder filtration of the contravariant defect. The outgoing case follows by duality. This gives a new proof of the implication from local finiteness to generation by Auslander--Reiten relations. Our argument uses defects and induction on their composition length. In this way, it avoids the radical filtration argument used by Zhu and Zhuang. We then sssume that $\C$ has enough projectives and injectives, and prove that $\C$ is locally finite if and only if its projective and injective subcategories are locally finite, the Euler functions associated with Auslander--Reiten relations have finite support, stable and costable morphisms detect ordinary morphisms away from the projective and injective boundary, and $\Rel(\C)=\AR(\C)$.
Our results recover the classical module and triangulated cases. They also explain why the two one-sided support conditions are equivalent in the triangulated setting, but may differ in a general extriangulated category.
\\[0.2cm]
\end{abstract}

\pagestyle{myheadings}
\markboth{\rightline{\scriptsize Panyue Zhou and Bin Zhu\hspace{2.3mm}}}
         {\leftline{\scriptsize Auslander--Reiten relations and local finiteness}}

\section{Introduction}

Let $R$ be a commutative artin ring and let $\Lambda$ be an Artin $R$-algebra. Butler \cite{Butler1981} proved that finite representation type implies that the defining relations of $K_0(\modcat\Lambda)$ are generated by almost split sequences. Auslander \cite{Auslander1984} later proved the converse.
A similar question arises for triangulated categories. Xiao and Zhu \cite[Theorem~2.1]{XiaoZhu2002} proved that local finiteness implies generation by Auslander--Reiten triangles. Over a field, the converse was proved by Zhou \cite[Theorem~3.14 and Remark~3.15]{Zhou2021}. Haugland \cite{Haugland2022} obtained related criteria involving a generator or a cogenerator.

Extriangulated categories were introduced by Nakaoka and Palu \cite{NakaokaPalu2019}. They include exact categories, triangulated categories, and extension-closed  subcategories of triangulated categories. Auslander--Reiten theory for extriangulated categories was later developed by Iyama, Nakaoka and Palu \cite{IyamaNakaokaPalu2024}. For an extriangulated category $\C$, $\Rel(\C)$ denotes the subgroup of $K_0^{\mathrm{sp}}(\C)$ generated by all $\E$-triangle relations, and $\AR(\C)$ denotes the subgroup generated by Auslander--Reiten $\E$-triangle relations.
Over a field, Zhu and Zhuang \cite[Theorems~3.12 and~4.4]{ZhuZhuang2021} proved that every locally finite extriangulated category has Auslander--Reiten $\E$-triangles and that $\Rel(\C)=\AR(\C)$.
We give a new proof of these results using defects from one-sided finiteness condition on
Hom supports.
Defects in extriangulated categories have been studied by Ogawa \cite{Ogawa2021} and Enomoto \cite{Enomoto2021}. Our argument is closely related to the functorial method used by Enomoto \cite{Enomoto2019} for exact categories. The main difference from the proof of Zhu and Zhuang is the induction process. Their proof passes to morphisms in higher powers of the radical and stops when a suitable radical power vanishes. Our proof instead removes one simple composition factor of the contravariant defect at each step. Thus the induction terminates when the composition length of the defect is exhausted.
Although Zhu and Zhuang \cite[Theorem~4.4]{ZhuZhuang2021} assume two-sided local finiteness, the contravariant part of their argument already suggests a one-sided phenomenon. Our proof makes this point explicit. It also works over a commutative artin ring and requires no separate finiteness assumption on the bifunctor $\E$.

 The converse problem is to determine when
$\Rel(\C)=\AR(\C)$ implies that $\C$ is locally finite.
The converse fails without further assumptions, see \cite[Remark~2]{MarcosMerklenPlatzeck2000} and \cite[Example~4.5(4)]{ZhuZhuang2021}.
To obtain a converse, two difficulties must be taken into account. For an Auslander--Reiten $\E$-triangle, its Euler relation may contain extra contributions when an inflation is not a monomorphism or a deflation is not an epimorphism. The finite Euler support condition (EFS) (see Definition \ref{def:EFS}) is introduced to control this phenomenon. There exists also a second difficulty. Stable and costable morphisms do not detect morphisms that factor through projective or injective objects. The conditions (${\rm D}_-$) and
(${\rm D}_+$) (see Definition \ref{def:detectors}) allow us to pass from stable and costable Hom spaces to ordinary Hom spaces.

We denote by $\Pcal=\Proj\C$ and $\Ical=\Inj\C$ the full subcategories of projective and injective objects, respectively. Recall that a full additive Krull--Schmidt subcategory is locally finite if, for every indecomposable object, only finitely many indecomposable objects in the subcategory admit a nonzero morphism to or from it, up to isomorphism (see Definition \ref{locall} for the precise definition).  It is useful to consider the two directions of local finiteness separately. The contravariant condition alone implies $\Rel(\C)=\AR(\C)$
by an argument using defects. The covariant statement follows by duality. These results are proved in Section~\ref{sec:one-sided}. Together with the detection conditions, they lead to the characterization of local finiteness proved in Section~\ref{sec:detectors}.

\begin{theorem}\label{thm:intro-main}{\rm (see Theorem \ref{thm:main} for details)}
Let $\C$ be a skeletally small Hom-finite Krull--Schmidt $R$-linear
extriangulated category with enough projectives and enough injectives.  Then
the following statements are equivalent.
\begin{enumerate}[label=\textup{(\arabic*)},leftmargin=2.7em]
\item $\C$ is locally finite;
\item Both $\Pcal$ and $\Ical$ are locally finite, $\C$ satisfies \emph{(EFS)}, \emph{(D$_-$)} and
\emph{(D$_+$)}, and
$\Rel(\C)=\AR(\C)$.
\end{enumerate}
\end{theorem}

For the implication \textup{(1)}$\Rightarrow$\textup{(2)}, the one-sided
defect theorem gives the equality of the relation subgroups.  The other
conditions follow directly from local finiteness.  For the converse, syzygy
and cosyzygy relations, together with $\Rel(\C)=\AR(\C)$ and (EFS), show that
stable and costable representable functors have finite support.  The
detection conditions then allow us to recover ordinary Hom spaces.  The
local finiteness of $\Pcal$ and $\Ical$ controls the projective and injective
boundary.  A useful sufficient condition for the detection assumptions is
that $\underline\C$ has a weak cogenerator and $\overline\C$ has a weak
generator.  By Lemma~\ref{lem:weak-detectors}, we give the precise statement.
The boundary assumptions are automatic when $\C$ has a projective generator
and an injective cogenerator.  In this case, $\Pcal$ and $\Ical$ are of
finite type by Corollary~\ref{cor:generator-cogenerator}.  The assumptions on the projective and injective boundaries are essential.  Indeed, the category in Example~\ref{ex:one-sided-not-self-dual} has enough projectives and injectives, satisfies (EFS) and both detection conditions, and has $\Rel(\C)=\AR(\C)=0$.  Nevertheless, the category is not locally finite, since its projective and injective boundaries fail to be locally finite.

The theorem recovers the classical module and triangulated cases.  For
$\modcat\Lambda$, the semisimple module $\Lambda/\rad\Lambda$ is a weak
cogenerator of the stable category.  The dual statement holds in the
costable category.  This is a special case of Enomoto
\cite[Proposition~4.7]{Enomoto2019}.  In a triangulated category, the
projective and injective objects vanish, so both detection conditions are
automatic. We also obtain a consequence for Frobenius extriangulated categories.  Let
$\C$ be a skeletally small Hom-finite Krull--Schmidt $R$-linear Frobenius
extriangulated category and assume that $\Rel(\C)=\AR(\C)$.  Then the same
equality holds in the triangulated stable category $\underline{\C}$, and
hence $\underline{\C}$ is locally finite.  For Frobenius exact categories,
Haugland \cite[Lemma~3.2]{Haugland2022} proved the corresponding passage to
the stable category.  By Proposition~\ref{prop:frobenius-stable-relations},
we extend this result to Frobenius extriangulated categories.

For $0$-Auslander extriangulated categories whose projective and injective
subcategories are of finite type, descriptions of the stable and costable
quotients by module categories make the detection conditions automatic.  In
particular, we obtain a converse criterion for the relative extriangulated
structure associated with a cluster tilting object in a Hom-finite
$2$-Calabi--Yau triangulated category.  We also give examples beyond the
exact and triangulated settings.  These include a locally finite category of
infinite type and a graded example whose projective and injective
subcategories have infinitely many indecomposable objects.

The paper is organized as follows.  In Section~\ref{sec:prelim}, we fix the notation and recall the required facts on extriangulated categories, Grothendieck groups and local finiteness.  In Section~\ref{sec:one-sided}, we use contravariant defects to prove that either one-sided local finiteness condition implies $\Rel(\C)=\AR(\C)$.  In Section~\ref{sec:euler}, we introduce the Euler functions and the finite support condition.  In Section~\ref{sec:detectors}, we introduce the detection conditions and prove Theorem~\ref{thm:intro-main}. In Section~~\ref{sec:applications}, we treat the module, triangulated and Frobenius cases. We also give examples of genuinely extriangulated categories that are neither exact nor triangulated.  In Section~\ref{sec:serre}, we establish a criterion under Serre duality.

\section{Preliminaries}\label{sec:prelim}

Throughout this paper, $R$ is a commutative artin ring and $\C$ is a
skeletally small, $R$-linear Krull--Schmidt category.  We say that $\C$ is
\emph{Hom-finite} if $\C(X,Y)$ has finite length as an $R$-module for all
$X,Y\in\C$, and write $\ell_R(-)$ for this length.  We write $\ind\C$ for a
fixed set of representatives of the isomorphism classes of indecomposable
objects.  All subcategories are full, additive and closed under isomorphisms
and direct summands.  The case over a field $k$ is obtained by taking
$R=k$.

\vspace{2mm}

\hspace{-5mm}{\bf 2.1 Extriangulated categories.}~
We use the terminology of \cite{NakaokaPalu2019}.  An extriangulated category
is a triple $(\C,\E,\mathfrak s)$ consisting of an additive category $\C$, an
additive bifunctor
\[
 \E\colon \C^{\mathrm{op}}\times\C\longrightarrow\mathrm{Ab},
\]
and an additive realization $\mathfrak s$ satisfying the axioms in
\cite[Definition~2.12]{NakaokaPalu2019}.  A realization of
$\delta\in\E(C,A)$ will be written
\begin{equation}\label{eq:Etriangle}
 A\xrightarrow{a}B\xrightarrow{b}C\overset{\delta}{\dashrightarrow}
\end{equation}
and called an $\E$-triangle.  Associated with \eqref{eq:Etriangle} are exact
sequences
\begin{align*}
 \C(U,A)&\xrightarrow{\C(U,a)}\C(U,B)
 \xrightarrow{\C(U,b)}\C(U,C)\longrightarrow\E(U,A);\\
 \C(C,U)&\xrightarrow{\C(b,U)}\C(B,U)
 \xrightarrow{\C(a,U)}\C(A,U)\longrightarrow\E(C,U).
\end{align*}
For morphisms $f\colon A\to A'$ and $g\colon C'\to C$, we use the standard notation
\[
 f_*\delta=\E(C,f)(\delta)~~\mbox{and}~~g^*\delta=\E(g,A)(\delta).
\]
An object $P$ is called projective if $\E(P,-)=0$, and an object $I$ is called injective if $\E(-,I)=0$.  We denote the full subcategories of projective and injective objects by $\Pcal=\Proj\C$ and $\Ical=\Inj\C$, respectively.
We say that an extriangulated category $\C$ \emph{has enough projectives} if, for every $X$, there is an $\E$-triangle
$$
 \Omega X\longrightarrow P_X\longrightarrow X\dashrightarrow,
 ~~\mbox{where}~P_X\in\Pcal.
$$
Dually, $\C$ \emph{has enough injectives} if, for every $X$, there is an $\E$-triangle
$$
 X\longrightarrow I^X\longrightarrow\Sigma X\dashrightarrow,
 ~~\mbox{where}~ I^X\in\Ical.
$$
The objects $\Omega X$ and $\Sigma X$ depend on these choices.
We do not require these assignments to be functorial.

Let $[\Pcal]$ be the ideal of morphisms which factor through a projective
object.  The stable category is $\underline\C=\C/[\Pcal]$.  Thus
\[
 \underline\C(X,Y)=\C(X,Y)/[\Pcal](X,Y).
\]
The costable category $\overline\C=\C/[\Ical]$ is defined dually.  Since
an ideal quotient has the same objects as the original category, we use the
same symbols for objects of $\C$ and their images in $\underline\C$ or
$\overline\C$.
\vspace{2mm}

\hspace{-5mm}{\bf 2.2 Auslander--Reiten \texorpdfstring{$\E$}{E}-triangles.}~
A morphism $b\colon B\to C$ is \emph{right almost split} if it is not a retraction
and every morphism $U\to C$ which is not a retraction factors through $b$.
Left almost split morphisms are defined dually.  An $\E$-triangle
\eqref{eq:Etriangle} is an \emph{Auslander--Reiten $\E$-triangle} if $a$ is
left almost split and $b$ is right almost split.  Its end terms are
indecomposable.  We refer to \cite{IyamaNakaokaPalu2024} for the general
Auslander--Reiten theory of extriangulated categories.

We shall use the following standard consequence of the Krull--Schmidt property.  If $U$ is indecomposable, then $\rad_\C(X,U)$ consists precisely of the morphisms $X\to U$ that are not retractions.  Dually, $\rad_\C(U,X)$ consists precisely of the morphisms $U\to X$ that are not sections.

From now on, unless otherwise stated, we assume that
$\C:=(\mathcal{C},\mathbb{E},\mathfrak{s})$ is an extriangulated category.

The next two elementary facts allow us to work in Section~\ref{sec:one-sided} without imposing any finiteness assumption on $\E$.

\begin{lemma}\label{lem:annihilation-AR}
Let $A,C\in\ind\C$ and let $0\ne\delta\in\E(C,A)$.  Suppose that
\[
 c^*\delta=0
 ~~\text{for every non-retraction }c\colon C'\to C.
\]
Then every realization
\[
 A\xrightarrow{x}B\xrightarrow{y}C\overset{\delta}{\dashrightarrow}
\]
is an Auslander--Reiten $\E$-triangle.
\end{lemma}

\begin{proof}
Since $\delta\ne0$, the morphism $y$ is not a retraction.  If
$c\colon C'\to C$ is not a retraction, then $c^*\delta=0$, and exactness of
\[
 \C(C',B)\xrightarrow{\C(C',y)}\C(C',C)
 \longrightarrow\E(C',A)
\]
shows that $c$ factors through $y$.  Thus $y$ is right almost split.

It remains to prove that $x$ is left almost split.  Let $a\colon A\to A'$ be a non-section.  Assume for a contradiction that $a_*\delta\ne0$, and choose a realization
$$
 A'\xrightarrow{x'}B'\xrightarrow{y'}C
 \overset{a_*\delta}{\dashrightarrow}.
$$
Then $y'$ is not a retraction, so the hypothesis gives
$(y')^*\delta=0$.  Exactness for the original $\E$-triangle yields a
morphism $s\colon B'\to B$ with $ys=y'$.  By $\rm (ET3)^{\rm op}$ there is a
morphism $d\colon A'\to A$ such that $(d,s,1_C)$ is a morphism of
$\E$-triangles from the last triangle to the original one.  Hence
\[
 \delta=d_*(a_*\delta)=(da)_*\delta.
\]
Because $a$ is not a section, $da$ is not invertible.  The ring
$\End_\C(A)$ is local, so $1_A-da$ is invertible.  On the other hand,
\[
 (1_A-da)_*\delta=\delta-(da)_*\delta=0,
\]
and applying the inverse of $1_A-da$ gives $\delta=0$, a contradiction.
Thus $a_*\delta=0$.  Exactness of
\[
 \C(B,A')\xrightarrow{\C(x,A')}\C(A,A')
 \longrightarrow\E(C,A')
\]
shows that every non-section $a\colon A\to A'$ factors through $x$.
Since $\delta\ne0$, $x$ is not a section.  Therefore $x$ is left almost
split, and the realization is an Auslander--Reiten $\E$-triangle.
\end{proof}

\hspace{-5mm}{\bf 2.3 Grothendieck groups and local finiteness.}~
The split Grothendieck group $K_0^{\mathrm{sp}}(\C)$ is the free abelian group
on $\ind\C$.  Equivalently, it is generated by symbols $[X]$, for $X\in \C$, subject to
$[X\oplus Y]=[X]+[Y]$, for $X,Y\in \C$.  To an $\E$-triangle
$\delta$ as in \eqref{eq:Etriangle} we associate its Euler relation
\[
 [\delta]=[A]+[C]-[B]\in K_0^{\mathrm{sp}}(\C).
\]
Let $\Rel(\C)$ be the subgroup generated by all such relations, and let
$\AR(\C)$ be the subgroup generated by the relations of Auslander--Reiten
$\E$-triangles.  Then
\[
 K_0(\C)=K_0^{\mathrm{sp}}(\C)/\Rel(\C),
\]
and the condition that the relations of $K_0(\C)$ are generated by
Auslander--Reiten $\E$-triangles is
\[
 \AR(\C)=\Rel(\C).
\]
We shall use the following comparison construction.  It is a direct
consequence of the extriangulated axiom (ET4).

\begin{lemma}\label{lem:comparison-triangle}
Let
\[
 \delta:\quad A\xrightarrow{a}B\xrightarrow{b}C
 \overset{\delta}{\dashrightarrow}
~~\mbox{and}~~
 \varepsilon:\quad A'\xrightarrow{x}E\xrightarrow{y}C
 \overset{\varepsilon}{\dashrightarrow}
\]
be $\E$-triangles.  Assume that $\varepsilon=\alpha_*\delta$ for some
$\alpha\colon A\to A'$.  Then there is a morphism
$\gamma\colon B\to E$ and an $\E$-triangle
\[
 \delta':\quad
 A\longrightarrow A'\oplus B
 \xrightarrow{(x,\gamma)}E\dashrightarrow
\]
such that $y\gamma=b$.  Moreover,
$
 [\delta]=[\varepsilon]+[\delta']
 ~~\text{in }K_0^{\mathrm{sp}}(\C).
$
\end{lemma}

\begin{proof}
Put
$
 \mu=\Delta_C^*(\delta\oplus\varepsilon)
 \in\E(C,A\oplus A'),
$
where $\Delta_C\colon C\to C\oplus C$ is the diagonal, and choose a
realization
$
 A\oplus A'\xrightarrow{j}M\xrightarrow{k}C
 \overset{\mu}{\dashrightarrow}.
$
Let $\iota_A,\iota_{A'}$ and $p_A,p_{A'}$ be the canonical inclusions and
projections.  Then
\[
 (p_A)_*\mu=\delta,
 ~~
 (p_{A'})_*\mu=\varepsilon.
\]
Apply (ET4) to the split $\E$-triangle
$A\xrightarrow{\iota_A}A\oplus A'\xrightarrow{p_{A'}}A'\dashrightarrow$
and the chosen realization of $\mu$.  We obtain a commutative diagram
\begin{equation}\label{diag:comparison-one}
\begin{array}{l}
\xymatrix@C=3.5em@R=2.6em{
 A \ar[r]^{\iota_A}\ar@{=}[d] &
 A\oplus A' \ar[r]^{p_{A'}}\ar[d]^{j} &
 A' \ar[d]^{x'} \\
 A \ar[r]_{m_1} &
 M \ar[r]_{e_1'}\ar[d]_{k} &
 E' \ar[d]^{y'} \\
 & C \ar@{=}[r] & C
}
\end{array}
\end{equation}
in which the first row is split, the middle row and the right column are
$\E$-triangles, and the middle column is the chosen realization of $\mu$.
The extension of the right column is $(p_{A'})_*\mu=\varepsilon$.
If $\theta_1\in\E(E',A)$ denotes the extension of the middle row, the
compatibility in (ET4) gives
\[
 (\iota_A)_*\theta_1=(y')^*\mu.
\]
Applying $(p_A)_*$ yields
$
 \theta_1=(y')^*(p_A)_*\mu=(y')^*\delta.
$
Since $A'\xrightarrow{x}E\xrightarrow{y}C$ realizes the same extension as
the right column, we may identify $E'$ with $E$.  After this identification
we write the middle morphism as $e_1\colon M\to E$.  Thus
\[
 A\xrightarrow{m_1}M\xrightarrow{e_1}E
 \overset{y^*\delta}{\dashrightarrow},
~~
 e_1j\iota_{A'}=x,
~~
 ye_1=k.
\]
Apply (ET4) a second time, now to the split $\E$-triangle
$A'\xrightarrow{\iota_{A'}}A\oplus A'\xrightarrow{p_A}A\dashrightarrow$
and the same realization of $\mu$.  This gives
\begin{equation}\label{diag:comparison-two}
\begin{array}{l}
\xymatrix@C=3.5em@R=2.6em{
 A' \ar[r]^{\iota_{A'}}\ar@{=}[d] &
 A\oplus A' \ar[r]^{p_A}\ar[d]^{j} &
 A \ar[d]^{a} \\
 A' \ar[r]_{m_2} &
 M \ar[r]_{e_2}\ar[d]_{k} &
 B \ar[d]^{b} \\
 & C \ar@{=}[r] & C
}
\end{array}
\end{equation}
where the middle row is an $\E$-triangle.  If its extension is
$\theta_2\in\E(B,A')$, the same compatibility gives
$
 (\iota_{A'})_*\theta_2=b^*\mu,
$
and hence
$
 \theta_2=b^*(p_{A'})_*\mu=b^*\varepsilon.
$
The two diagrams show, in particular,
$
 e_1m_2=x,
 e_2m_1=a,
 ye_1=be_2.
$
Since $\varepsilon=\alpha_*\delta$, functoriality of $\E$ gives
$
 b^*\varepsilon
 =b^*(\alpha_*\delta)
 =\alpha_*(b^*\delta)=0.
$
Thus the middle row of \eqref{diag:comparison-two} is split.  Choose an isomorphism
$\varphi\colon M\xrightarrow{\sim}A'\oplus B$ such that
\[
 \varphi m_2=\binom{1}{0},
~
 e_2=(0,1)\varphi.
\]
Write
$
 e_1\varphi^{-1}=(x,\gamma)
$
for some $\gamma\colon B\to E$; its first component is $x$ because
$e_1m_2=x$.  The identity $ye_1=be_2$ now gives $y\gamma=b$.  Transporting
the middle row of \eqref{diag:comparison-one} along $\varphi$ gives the required $\E$-triangle
$\delta'$.  Finally,
\[
 [\varepsilon]+[\delta']
 =([A']+[C]-[E])+([A]+[E]-[A']-[B])
 =[A]+[C]-[B]=[\delta].
\]
\end{proof}

For $X\in\C$, put
\begin{align*}
 \Supp\C(-,X)&=\{U\in\ind\C\mid\C(U,X)\ne0\},\\
 \Supp\C(X,-)&=\{U\in\ind\C\mid\C(X,U)\ne0\}.
\end{align*}

Recall the notion of local finiteness from \cite[Definition 2.1]{XiaoZhu2002} and
\cite[Definition~3.1]{ZhuZhuang2021}.  We use the same terminology for full
additive Krull--Schmidt subcategories, with the Hom spaces inherited from
$\C$.

\begin{definition} \label{locall}
The category $\C$ is \emph{of finite type}, if ind$\C$ is a finite set.
The category $\C$ is \emph{locally finite} if, for every indecomposable
$X\in\C$, both $\Supp\C(-,X)$ and $\Supp\C(X,-)$ are finite.
More generally, if $\A$ is a full additive Krull--Schmidt subcategory of
$\C$, we say that $\A$ is \emph{locally finite} if, for every
$X\in\ind\A$, both sets
\[
 \{U\in\ind\A\mid\C(U,X)\ne0\}
 \quad\text{and}\quad
 \{U\in\ind\A\mid\C(X,U)\ne0\}
\]
are finite.  Equivalently, $\A$ is locally finite as a category with the
Hom spaces inherited from $\C$.
\end{definition}

For later use, we separate the two sides of local finiteness.  We write
\begin{align*}
 (\mathrm{LF}_-)&:\quad |\Supp\C(-,X)|<\infty
 \quad\text{for every }X\in\ind\C,\\
 (\mathrm{LF}_+)&:\quad |\Supp\C(X,-)|<\infty
 \quad\text{for every }X\in\ind\C.
\end{align*}
Thus $\C$ is locally finite if and only if both $(\mathrm{LF}_-)$ and
$(\mathrm{LF}_+)$ hold.

Local finiteness is weaker than finite type.  A concrete example is provided by the derived category $D^b(kQ)$ of a Dynkin
quiver $Q$ over an algebraically closed field $k$. It is locally finite, but  it is not of finite type.

\section{One-sided local finiteness and Auslander--Reiten relations}\label{sec:one-sided}

No assumption that $\C$ has enough projectives or injectives is needed in this section.  We now give the proof using defects announced in the introduction.  Contravariant defects in extriangulated categories have been studied by Ogawa \cite{Ogawa2021} and Enomoto \cite{Enomoto2021}.  We only use their elementary presentation and prove the reduction argument directly.  Zhu and Zhuang \cite[Theorem~4.4]{ZhuZhuang2021} reduce relations by passing to morphisms in successively higher powers of the radical.  Their argument terminates when a suitable radical power vanishes.  Our proof uses a composition series of the contravariant defect and removes one simple factor at each step.  This argument is closely related to the functorial approach used by Enomoto \cite{Enomoto2019} for exact categories.

For an $\E$-triangle
$$
 \delta:\quad A\xrightarrow{a}B\xrightarrow{b}C\dashrightarrow,
$$
we define its contravariant defect as the additive contravariant functor
$$
 \widetilde\delta
 :=\coker\bigl(\C(-,B)\xrightarrow{\C(-,b)}\C(-,C)\bigr).
$$
For $U\in\ind\C$, set
$$
 S_U:=\C(-,U)/\rad_\C(-,U).
$$
Since $\C$ is Krull--Schmidt, the functor $S_U$ is simple.  We shall use the radical criterion from Section~\ref{sec:prelim} throughout this section.

\begin{lemma}\label{lem:one-sided-AR}
Assume $\C$ satisfies $(\mathrm{LF}_-)$.  Then every nonprojective object
$U\in\ind\C$ is the end term of an Auslander--Reiten $\E$-triangle.
\end{lemma}

\begin{proof}
Fix a nonprojective $U\in\ind\C$.  By $(\mathrm{LF}_-)$, there are only
finitely many indecomposable objects $X_1,\cdots,X_t$ with
$\C(X_i,U)\ne0$.  Since $\C$ is Hom-finite, each $\rad_\C(X_i,U)$ has finite length as an
$R$-module and is therefore finitely generated.
Choose generators $f_{ij}$ and let
$$
 q\colon~ Q:=\bigoplus_{i,j}X_i\longrightarrow U
$$
be the morphism whose components are the $f_{ij}$.  Let $h:X\to U$ be a
non-retraction and decompose $X$ into indecomposable summands.  No component of $h$
is a retraction.  After identifying each nonzero indecomposable summand with
one of the chosen representatives $X_i$, every component lies in
$\rad_\C(X_i,U)$ and is an $R$-linear combination of the chosen generators.
Hence $h$ factors through $q$.  Thus every non-retraction to $U$ factors
through $q$.

Since $U$ is nonprojective, there exists an object $A_0\in\C$ such that
$
\E(U,A_0)\neq0.
$
Choose a nonzero extension $\delta_0\in\E(U,A_0)$ and a realization
\[
A_0\longrightarrow B_0\xrightarrow{~p~}U
\overset{\delta_0}{\dashrightarrow}.
\]
Since $\delta_0\neq0$, this $\E$-triangle is nonsplit.
The morphism
$
 (p,~q):B_0\oplus Q\longrightarrow U
$
is a deflation.  Indeed, it is the composite
\[
 B_0\oplus Q\xrightarrow{p\oplus1_Q}U\oplus Q
 \xrightarrow{(1_U,~q)}U,
\]
where the first morphism is a deflation and the second is a split deflation.
Deflations are closed under composition by \cite[Remark~2.16]{NakaokaPalu2019}.
Since the $\E$-triangle is nonsplit,
$p$ is not a retraction, and hence $p\in\rad_\C(B_0,U)$.  By construction,
$q\in\rad_\C(Q,U)$.  Thus $(p,q)$ belongs to
$\rad_\C(B_0\oplus Q,U)$ and is not a retraction.  Since every
non-retraction to $U$ factors through $q$, it also factors through $(p,q)$.
Hence $(p,q)$ is a right almost split deflation.

Let
\[
 K\longrightarrow B_0\oplus Q\xrightarrow{~(p,~q)~}U
 \overset{\eta}{\dashrightarrow}
\]
be the resulting nonsplit $\E$-triangle.  Then $c^*\eta=0$ for every
non-retraction $c:U'\to U$.  Write $K=K_1\oplus\cdots\oplus K_s$ with the
$K_i$ indecomposable and choose a component $\eta_i\in\E(U,K_i)$ which is
nonzero.  For every non-retraction $c\colon U'\to U$ one has
$c^*\eta_i=0$, because $\eta_i$ is a component of $\eta$.  Since
$K_i$ and $U$ are indecomposable and $\eta_i\ne0$,
it follows from Lemma~\ref{lem:annihilation-AR} that  any realization of
$\eta_i$ is an Auslander--Reiten $\E$-triangle ending at $U$.
\end{proof}

\begin{lemma}[Defect reduction]\label{lem:defect-reduction}
Assume $\C$ satisfies $(\mathrm{LF}_-)$.  Let
\[
 \delta:\quad A\xrightarrow{a}B\xrightarrow{b}C\dashrightarrow
\]
be a nonsplit $\E$-triangle, and let
$\widetilde\delta\twoheadrightarrow S_U$ be a simple quotient.  Then $U$ is
nonprojective.  Let
\[
 \rho_U:\quad A_U\xrightarrow{x_U}E_U\xrightarrow{y_U}U\dashrightarrow
\]
be an Auslander--Reiten $\E$-triangle ending at $U$, which exists by
{\rm Lemma~\ref{lem:one-sided-AR}}.  Then there exists an $\E$-triangle $\delta'$ such that
\[
 [\delta]=[\rho_U]+[\delta']
\]
and there is a short exact sequence
\[
 0\longrightarrow\widetilde{\delta'}\longrightarrow
 \widetilde\delta\longrightarrow S_U\longrightarrow0.
\]
\end{lemma}

\begin{proof}
Let $\pi:\widetilde\delta\twoheadrightarrow S_U$ be the chosen quotient.
The composite
\[
 \C(-,C)\longrightarrow\widetilde\delta\xrightarrow{\pi}S_U
\]
is nonzero.  By Yoneda, it is represented by a morphism $p:C\to U$ whose
class in
\[
 S_U(C)=\C(C,U)/\rad_\C(C,U)
\]
is nonzero.  Thus $p\notin\rad_\C(C,U)$. By the radical criterion from Section~\ref{sec:prelim}, the morphism $p$ is a retraction.  Choose a section
$i:U\to C$ with $pi=1_U$.  Since $pi=1_U$, the endomorphism $ip$ of $C$ is idempotent. As $\C$ is Krull--Schmidt, this idempotent splits. Hence $C\simeq U\oplus C'$ for some object $C'$.
After transporting $\delta$ and the quotient $\pi$ along this isomorphism,
we may write
\[
 C=U\oplus C',~~ p=(1_U,0),~~
 b=\binom{b_U}{b'}.
\]
Let $q_\delta\colon\C(-,C)\twoheadrightarrow\widetilde\delta$ be the canonical
quotient, and let $\theta_p:\C(-,C)\to S_U$ be the natural transformation
corresponding under Yoneda to the class $[p]\in S_U(C)$.  Thus
$\theta_{p,X}(c)=[pc]$ for $c:X\to C$.  By the choice of $p$, we have
$\pi q_\delta=\theta_p$.  Evaluating this equality at $B$ on the element
$b\in\C(B,C)$ gives
$
 0=(\pi q_\delta)_B(b)=[pb]=[b_U]\in S_U(B).
$
Thus $b_U\in\rad_\C(B,U)$. In particular, $b_U$ is not a retraction.  Moreover,
we also claim that $U$ is nonprojective. Indeed, if $U$ were projective, exactness of
$\C(U,B)\to\C(U,C)\to\E(U,A)$ would imply
$\widetilde\delta(U)=0$, contradicting
$\widetilde\delta(U)\twoheadrightarrow S_U(U)\ne0$.

Since $y_U$ is right almost split, there is $\beta_U:B\to E_U$ with
$b_U=y_U\beta_U$.  Take the direct sum of $\rho_U$ with the split
$\E$-triangle on $C'$ and write it as
\[
 \varepsilon:\quad
 A_U\xrightarrow{\bar x_U}E:=E_U\oplus C'
 \xrightarrow{y:=y_U\oplus1_{C'}}C=U\oplus C'\dashrightarrow,
 ~~ \bar x_U=\binom{x_U}{0}.
\]
Put
\[
 \beta=\binom{\beta_U}{b'}:B\longrightarrow E.
\]
Then $b=y\beta$.  By (ET3)$^{\mathrm{op}}$, there is a morphism
$\alpha:A\to A_U$ such that $(\alpha,\beta,1_C)$ is a morphism from
$\delta$ to $\varepsilon$.  Thus $\varepsilon=\alpha_*\delta$. By Lemma~\ref{lem:comparison-triangle},
we give a morphism $\gamma:B\to E$ with $y\gamma=b$ and an $\E$-triangle
\[
 \delta':\quad
 A\longrightarrow A_U\oplus B
 \xrightarrow{(\bar x_U,\gamma)}E\dashrightarrow.
\]
Since the split summand contributes no relation, $[\varepsilon]=[\rho_U]$. It follows again from Lemma~\ref{lem:comparison-triangle} that
\[
 [\delta]=[\rho_U]+[\delta'].
\]

For $X\in\C$, composition with $y$ induces
\[
 \Phi_X:\widetilde{\delta'}(X)\longrightarrow\widetilde\delta(X),
 ~~ [h]\longmapsto[yh].
\]
It is well defined because
$y(\bar x_U,\gamma)(u,v)=bv$.  If $\Phi_X([h])=0$, then
$yh=bv=y\gamma v$ for some $v:X\to B$.  Hence
$y(h-\gamma v)=0$.  Exactness of the Hom sequence for $\varepsilon$ gives
$h-\gamma v=\bar x_Uu$ for some $u:X\to A_U$, so $[h]=0$.  Thus
$\Phi_X$ is injective.

The projection $p:C=U\oplus C'\to U$ induces a natural epimorphism
\[
 \Psi:\widetilde\delta\longrightarrow
 \C(-,U)/\operatorname{Im}\C(-,y_U).
\]
Indeed, at an object $X$ it is given by
$
 \Psi_X([c])=[pc].
$
It is well defined because $pb=b_U=y_U\beta_U$, and it is surjective
because $p$ is a retraction.  If $[c]$ lies in its kernel, write
$c=(c_U,c')$ and choose $e_U:X\to E_U$ with $c_U=y_Ue_U$.  Then
$h=(e_U,c'):X\to E$ satisfies $yh=c$, so $[c]$ lies in the image of
$\Phi_X$.  The reverse inclusion is immediate.  Therefore
\[
 0\longrightarrow\widetilde{\delta'}
 \longrightarrow\widetilde\delta
 \xrightarrow{\Psi}
 \C(-,U)/\operatorname{Im}\C(-,y_U)\longrightarrow0
\]
is exact.  Since $y_U$ is right almost split,
$
 \operatorname{Im}\C(-,y_U)=\rad_\C(-,U),
$
so the last term is $S_U$.  Under this identification, $\Psi$ is precisely
the chosen quotient $\pi$.  In fact, if
$$q_\delta:\C(-,C)\twoheadrightarrow\widetilde\delta$$
denotes the
canonical quotient, then both composites in
\[
\xymatrix@C=4.5em{
 \C(-,C) \ar@{->>}[r]^{q_\delta} \ar[dr]_{\theta_p} &
 \widetilde\delta \ar@{->>}[d]^{\pi} \\
 & S_U
}
\]
are equal to $\theta_p$.  Since $q_\delta$ is an epimorphism,
$\Psi=\pi$, and the asserted short exact sequence follows.
\end{proof}

\begin{theorem}\label{thm:one-sided-relations}
Let $\C$ be a skeletally small Hom-finite Krull--Schmidt $R$-linear
extriangulated category.  If $\C$ satisfies $(\mathrm{LF}_-)$ or $(\mathrm{LF}_+)$, then
\[
 \Rel(\C)=\AR(\C).
\]
\end{theorem}

\begin{proof}
Assume that $(\mathrm{LF}_-)$ holds.  The case $(\mathrm{LF}_+)$ follows by applying
the same argument to the opposite extriangulated category
$\C^{\mathrm{op}}$.  An $\E$-triangle $$A\to B\to C\dashrightarrow$$ in $\C$ becomes an $\E$-triangle
$C\to B\to A\dashrightarrow$ in $\C^{\mathrm{op}}$, so its relation
$[A]+[C]-[B]$ is unchanged, and Auslander--Reiten $\E$-triangles correspond under
opposition.  Hence, under the natural identification
$K_0^{\mathrm{sp}}(\C^{\mathrm{op}})=K_0^{\mathrm{sp}}(\C)$, the relation
and Auslander--Reiten subgroups correspond to $\Rel(\C)$ and $\AR(\C)$,
respectively.
Let
\[
 \delta:\quad A\longrightarrow B\longrightarrow C\dashrightarrow
\]
be an $\E$-triangle.  Since $\Supp\C(-,C)$ is finite and every $\C(X,C)$ has finite
$R$-length, the representable functor $\C(-,C)$ has finite length.  Indeed,
for a subfunctor $F\subseteq\C(-,C)$ set
\[
 \nu(F)=\sum_{X\in\Supp\C(-,C)}\ell_R F(X).
\]
A strict inclusion of subfunctors strictly increases $\nu$.  Indeed, if
$F\subsetneq G$, then $F(Y)\subsetneq G(Y)$ for some object $Y$.  Decompose
$Y$ into indecomposable summands.  By additivity,
$F(X)\subsetneq G(X)$ for at least one indecomposable summand $X$.  Since
$G(X)\ne0$, such an $X$ belongs to $\Supp\C(-,C)$.
 Dually, a strict descending
chain strictly decreases $\nu$.  Since
$0\le\nu(F)\le\nu(\C(-,C))$, neither type of chain can be infinite.
Hence $\C(-,C)$, and therefore its quotient $\widetilde\delta$, has finite
length.  Write $\length(F)$ for
the composition length of a functor $F$ of finite length.

We argue by induction on $\length(\widetilde\delta)$.  If the length is
zero, then $\widetilde\delta=0$.  Evaluating at $C$ shows that $1_C$
factors through the deflation $B\to C$, so $\delta$ is split and
$[\delta]=0$.

Suppose $\length(\widetilde\delta)>0$.  Choose a simple quotient $F$ of
$\widetilde\delta$.  We claim that $F\simeq S_U$ for some $U\in\ind\C$.
Choose an object on which $F$ is nonzero and decompose it into indecomposable
summands.  Since $F$ is additive, $F(U)\ne0$ for some indecomposable summand
$U$.  Choose a nonzero element $x\in F(U)$.
 By the Yoneda lemma, $x$
corresponds to a natural transformation
\[
\varphi:\C(-,U)\longrightarrow F,
~~
\varphi_X(f)=F(f)(x)
\]
for $f:X\to U$.  Since
$
\varphi_U(1_U)=x\neq0,
$
the natural transformation $\varphi$ is nonzero.  As $F$ is simple,
its image is a nonzero subfunctor of $F$, and hence
$\operatorname{Im}\varphi=F$.  Thus $\varphi$ is an epimorphism.
Let $G$ be its kernel.  If $g\in G(X)\subseteq\C(X,U)$ were a retraction,
choose $s:U\to X$ with $gs=1_U$.  Naturality would then give
$1_U\in G(U)$, forcing $G=\C(-,U)$, a contradiction.  Hence
$G\subseteq\rad_\C(-,U)$.  Since $F$ is simple, $G$ is a maximal proper
subfunctor of $\C(-,U)$.  The radical is proper because
$1_U\notin\rad_\C(U,U)$, and therefore maximality forces
$G=\rad_\C(-,U)$.  Thus $F\simeq S_U$.
By Lemma~\ref{lem:defect-reduction}, we give
\[
 [\delta]=[\rho_U]+[\delta']
\]
and
\[
 \length(\widetilde{\delta'})=\length(\widetilde\delta)-1.
\]
By induction, $[\delta']\in\AR(\C)$, hence $[\delta]\in\AR(\C)$.
Thus $\Rel(\C)\subseteq\AR(\C)$, and the reverse inclusion is immediate.
\end{proof}

\begin{example}\label{ex:one-sided-not-self-dual}
The conditions $(\mathrm{LF}_-)$ and $(\mathrm{LF}_+)$ are not equivalent
for general extriangulated categories, even in the presence of enough
projectives and injectives.  Let $k$ be a field and let $\C$ be the
Krull--Schmidt $k$-linear additive category with indecomposable objects
$X_1,X_2,\cdots$ and
\[
 \C(X_i,X_j)=
 \begin{cases}
 k,& i\le j,\\
 0,& i>j,
 \end{cases}
\]
with the evident composition.  Equip $\C$ with the split extriangulated
structure.  Then every object is both projective and injective.  For each
$j$,
\[
 \Supp\C(-,X_j)=\{X_1,\cdots,X_j\}
\]
is finite, whereas
\[
 \Supp\C(X_i,-)=\{X_i,X_{i+1},X_{i+2},\cdots\}
\]
is infinite.  Thus $(\mathrm{LF}_-)$ holds but $(\mathrm{LF}_+)$ does not.
Passing to the opposite category shows that the converse implication also
fails.
\end{example}


\section{Euler functions}\label{sec:euler}

For $U\in\C$, define homomorphisms
\[
 \chi_U^-,\chi_U^+\colon K_0^{\mathrm{sp}}(\C)\longrightarrow\mathbb Z
\]
by
\[
 \chi_U^-([X])=\ell_R\bigl(\C(U,X)\bigr),~~
 \chi_U^+([X])=\ell_R\bigl(\C(X,U)\bigr).
\]
For $r\in K_0^{\mathrm{sp}}(\C)$, write
\begin{align*}
 \Supp\chi^-(r)&=\{U\in\ind\C\mid\chi_U^-(r)\ne0\},\\
 \Supp\chi^+(r)&=\{U\in\ind\C\mid\chi_U^+(r)\ne0\}.
\end{align*}

The following elementary identity will be used repeatedly.

\begin{lemma}\label{lem:euler-defect}
Let $A\xrightarrow{a}B\xrightarrow{b}C\dashrightarrow$ be an
$\E$-triangle and let $U\in\C$.  Then
\begin{align*}
 \chi_U^-([A]+[C]-[B])
 &=\ell_R\bigl(\ker\C(U,a)\bigr)
   +\ell_R\bigl(\coker\C(U,b)\bigr),\\
 \chi_U^+([A]+[C]-[B])
 &=\ell_R\bigl(\ker\C(b,U)\bigr)
   +\ell_R\bigl(\coker\C(a,U)\bigr).
\end{align*}
In particular, both integers are nonnegative.
\end{lemma}

\begin{proof}
Set $V_1=\C(U,A)$, $V_2=\C(U,B)$ and $V_3=\C(U,C)$.  Exactness at
$V_2$ of the first Hom sequence displayed after \eqref{eq:Etriangle} gives
\[
 \ell_R(V_1)-\ell_R(V_2)+\ell_R(V_3)
 =\ell_R\bigl(\ker(V_1\to V_2)\bigr)
  +\ell_R\bigl(\coker(V_2\to V_3)\bigr).
\]
This proves the first identity.  The second is proved in the same way, using
exactness of the second Hom sequence associated with \eqref{eq:Etriangle}.
\end{proof}

\begin{definition}\label{def:EFS}
We say that $\C$ satisfies {(EFS)} if, for every Auslander--Reiten
$\E$-triangle $\delta$, both $\Supp\chi^-([\delta])$ and
$\Supp\chi^+([\delta])$ are finite.
\end{definition}

The property is automatic for exact categories.

\begin{proposition}\label{prop:exact-EFS}
Every Hom-finite Krull--Schmidt exact category has \emph{(EFS)}.  More
precisely, if
\[
 \delta:\quad 0\longrightarrow A\xrightarrow{~a~}B
 \xrightarrow{~b~}C\longrightarrow0
\]
is an Auslander--Reiten conflation, then
\[
 \Supp\chi^-([\delta])=\{C\}~~\mbox{and}~~
 \Supp\chi^+([\delta])=\{A\}.
\]
Here the sets are understood as sets of isomorphism classes.
\end{proposition}

\begin{proof}
Since $a$ is a monomorphism, $\C(U,a)$ is injective.
By Lemma~\ref{lem:euler-defect}, we give
\[
 \chi_U^-([\delta])=\ell_R\bigl(\coker\C(U,b)\bigr).
\]
If $U$ is indecomposable and $U\not\simeq C$, no morphism $U\to C$ is a
retraction.  Since $b$ is right almost split, every such morphism factors
through $b$, and the cokernel vanishes.  For $U=C$, the identity of $C$ does
not factor through $b$, since $b$ is not a retraction.  Hence the cokernel is
nonzero.  The second assertion is dual, using that $b$ is an epimorphism and
$a$ is left almost split.
\end{proof}

We also use the following criterion based on Serre duality.

\begin{proposition}\label{prop:Serre-EFS}
Let $E_R$ be an injective cogenerator of $\modcat R$ and put
$D={\rm Hom}_R(-,E_R)$.  Assume that $\C$ admits an $R$-linear Serre functor $S$,
that is, an $R$-linear autoequivalence together with bifunctorial
isomorphisms
\[
 D\C(X,Y)\simeq\C(Y,SX).
\]
Recall that $D$ is an exact duality on $R$-modules of finite length. In
particular, it preserves vanishing and composition length.
For an Auslander--Reiten $\E$-triangle
$\delta\colon A\xrightarrow{a}B\xrightarrow{b}C\dashrightarrow$, one has
\begin{align*}
 \Supp\chi^-([\delta])&\subseteq\{C,S^{-1}A\},\\
 \Supp\chi^+([\delta])&\subseteq\{A,SC\}.
\end{align*}
Consequently, $\C$ has \emph{(EFS)}.
\end{proposition}

\begin{proof}
By Lemma~\ref{lem:euler-defect}, $\chi_U^-([\delta])$ is the sum of a
kernel term and a cokernel term.  The latter is
$\ell_R(\coker\C(U,b))$, and it can be nonzero only for $U\simeq C$,
because $b$ is right almost split.

By Serre duality, the dual of
$$
 \C(U,A)\xrightarrow{\C(U,a)}\C(U,B)
$$
is identified with
$$
 \C(B,SU)\xrightarrow{\C(a,SU)}\C(A,SU).
$$
It follows that $\ker\C(U,a)$ is nonzero precisely when the latter map has
nonzero cokernel.  Since $a$ is left almost split, this can happen only when
$SU\simeq A$, or equivalently $U\simeq S^{-1}A$.  This proves the first
inclusion.  The second statement follows by the dual argument.  The cokernel of
$\C(a,U)$ is supported at $A$.  After Serre duality, the kernel of
$\C(b,U)$ is supported at $SC$.
\end{proof}

For triangulated categories one can be more explicit.  If
$A\to B\to C\to A[1]$ is an Auslander--Reiten triangle, the familiar
exact Hom sequences give
\begin{equation}\label{eq:tri-support}
 \Supp\chi^-([\delta])\subseteq\{C,C[-1]\}
~~\mbox{and}~~
 \Supp\chi^+([\delta])\subseteq\{A,A[1]\}.
\end{equation}
See, for example, \cite[Lemmas~2.2 and~2.3]{Haugland2022} or
\cite[Lemma~3.11]{Zhou2021}.  Thus a Hom-finite Krull--Schmidt triangulated
category has (EFS) whenever the Auslander--Reiten triangles under
consideration exist.

\section{Detection conditions}\label{sec:detectors}

Finite support of the Euler functions does not by itself control the
support of representable functors.  One standard way to detect
nonprojective objects in the stable category is to use a weak cogenerator.
An object $Q$ of $\underline\C$ is a \emph{weak cogenerator} if
$\underline\C(U,Q)=0$ implies that $U$ is projective.  Dually, an object
$T$ of $\overline\C$ is a \emph{weak generator} if
$\overline\C(T,U)=0$ implies that $U$ is injective.  For the converse
theorem, we only need these detection properties objectwise.  In this
form, the detecting object may depend on the support under consideration.

\begin{definition}\label{def:detectors}
We say that $\C$ satisfies (D$_-$) if, for every $X\in\ind\C$, there is an object
$Q_X\in\C$ such that, for every $U\in\ind\C$,
\[
 \C(U,X)\ne0
 ~~\Longrightarrow~~
 U\in\Pcal\ \text{or}\ \underline\C(U,Q_X)\ne0.
\]
Dually, $\C$ satisfies (D$_+$) if, for every $X\in\ind\C$, there is
an object $T_X\in\C$ such that, for every $U\in\ind\C$,
\[
 \C(X,U)\ne0
~~\Longrightarrow~~
 U\in\Ical\ \text{or}\ \overline\C(T_X,U)\ne0.
\]
\end{definition}

Thus a single weak cogenerator or weak generator gives a uniform form of the
corresponding detection condition.

\begin{lemma}\label{lem:weak-detectors}
If $\underline\C$ has a weak cogenerator, then $\C$ satisfies
\emph{(D$_-$)}.  If $\overline\C$ has a weak generator, then $\C$ satisfies
\emph{(D$_+$)}.
\end{lemma}

\begin{proof}
Let $Q$ be a weak cogenerator of $\underline\C$.  If
$\C(U,X)\ne0$ and $U$ is not projective, then
$\underline\C(U,Q)\ne0$ by definition.  Thus one may take $Q_X=Q$ for
all $X$, which is precisely {(D$_-$)}.  The second assertion is dual.
\end{proof}

The objects $Q_X$ and $T_X$ need not be chosen from a finite set.  This is
important in the triangulated setting, where local finiteness is genuinely
weaker than finite type.  In a module category, by contrast, a single object
works for every $X$.

In the next lemma, we relate these conditions to Grothendieck relations.

\begin{lemma}\label{lem:stable-defect}
Assume that $\C$ has enough projectives.  For $Q\in\C$, choose an
$\E$-triangle
\[
 \Omega Q\xrightarrow{u}P_Q\xrightarrow{p}Q\dashrightarrow,
 ~~\mbox{where}~P_Q\in\Pcal,
\]
and set
\[
 r_Q^-=[\Omega Q]+[Q]-[P_Q]\in\Rel(\C).
\]
Then, for every $U\in\C$,
\[
 \chi_U^-(r_Q^-)
 =\ell_R\bigl(\ker\C(U,u)\bigr)
  +\ell_R\bigl(\underline\C(U,Q)\bigr).
\]
In particular,
\[
 \underline\C(U,Q)\ne0~~\Longrightarrow~~
 \chi_U^-(r_Q^-)>0.
\]
Dually, if
\[
 T\xrightarrow{i}I^T\xrightarrow{v}\Sigma T\dashrightarrow,
~~\mbox{where}~I^T\in\Ical,
\]
and $r_T^+=[T]+[\Sigma T]-[I^T]$, then
\[
 \chi_U^+(r_T^+)
 =\ell_R\bigl(\ker\C(v,U)\bigr)
  +\ell_R\bigl(\overline\C(T,U)\bigr).
\]
\end{lemma}

\begin{proof}
By Lemma~\ref{lem:euler-defect}, we give
\[
 \chi_U^-(r_Q^-)
 =\ell_R\bigl(\ker\C(U,u)\bigr)
  +\ell_R\bigl(\coker\C(U,p)\bigr).
\]
We claim that the last cokernel is $\underline\C(U,Q)$.  Every morphism
in the image of $\C(U,p)$ factors through a projective object.
Conversely, suppose that $f\colon U\to Q$ factors as
$U\to P'\xrightarrow{g}Q$ with $P'$ projective.  Exactness of
\[
 \C(P',P_Q)\longrightarrow\C(P',Q)\longrightarrow\E(P',\Omega Q)=0
\]
shows that $g$ factors through $p$.  Hence $f$ lies in the image of
$\C(U,p)$.  This proves the first formula.

For the dual statement, the image of
$\C(I^T,U)\to\C(T,U)$ consists exactly of the morphisms $T\to U$ which
factor through an injective object.  The second formula follows from the dual argument and the second identity
in Lemma~\ref{lem:euler-defect}.
\end{proof}

\begin{lemma}\label{lem:locally-finite-detection}
Let $\C$ be a locally finite Hom-finite Krull--Schmidt extriangulated category
with enough projectives and injectives.  Then $\Pcal$ and $\Ical$ are locally
finite, and $\C$ satisfies \emph{(EFS)}, \emph{(D$_-$)} and \emph{(D$_+$)}.
\end{lemma}

\begin{proof}
The first assertion is immediate because $\Pcal$ and $\Ical$ are full
subcategories of $\C$. For an Auslander--Reiten $\E$-triangle
$A\to B\to C\dashrightarrow$, it follows from Lemma~\ref{lem:euler-defect} that
$$
 \Supp\chi^-([\delta])
 \subseteq\Supp\C(-,A)\cup\Supp\C(-,B)\cup\Supp\C(-,C),
$$
and the dual statement holds for $\chi^+$.
 Hence (EFS) follows from local finiteness.

Fix $X\in\ind\C$ and let $\mathcal S_X^-$ be the finite set of
nonprojective indecomposables $U$ with $\C(U,X)\ne0$.  Put
$Q_X=\bigoplus_{U\in\mathcal S_X^-}U$, with $Q_X=0$ if the set is empty.
For $U\in\mathcal S_X^-$, the canonical inclusion $U\to Q_X$ is nonzero in
$\underline\C$; otherwise $1_U$ would factor through a projective object,
forcing $U$ itself to be projective.  Thus (D$_-$) holds.  The proof of
(D$_+$) is dual.
\end{proof}

We now prove the main result of the paper.

\begin{theorem}\label{thm:main}
Let $\C$ be a skeletally small Hom-finite Krull--Schmidt $R$-linear
extriangulated category with enough projectives and enough injectives.  Then
the following conditions are equivalent.
\begin{enumerate}[label=\textup{(\arabic*)},leftmargin=2.7em]
\item $\C$ is locally finite.
\item Both $\Pcal$ and $\Ical$ are locally finite, $\C$ satisfies \emph{(EFS)}, \emph{(D$_-$)} and
\emph{(D$_+$)}, and
$$\Rel(\C)=\AR(\C).$$
\end{enumerate}
\end{theorem}

\begin{proof}
Assume first that $\C$ is locally finite.  By
Lemma~\ref{lem:locally-finite-detection}, $\Pcal$ and $\Ical$ are locally
finite and $\C$ satisfies (EFS), (D$_-$) and (D$_+$).  Moreover,
$(\mathrm{LF}_-)$ holds, so Theorem~\ref{thm:one-sided-relations} gives
$\Rel(\C)=\AR(\C)$.  Thus \textup{(1)} implies \textup{(2)}.

Conversely, assume \textup{(2)} and fix $X\in\ind\C$.  Choose an
$\E$-triangle
\[
 \Omega X\longrightarrow P^X\xrightarrow{p_X}X\dashrightarrow,
 ~~\mbox{where}~P^X\in\Pcal,
\]
and write $P^X=P_1\oplus\cdots\oplus P_t$ with $P_i\in\ind\Pcal$.
If $P\in\ind\Pcal$ and $0\ne f:P\to X$, projectivity gives
$g:P\to P^X$ with $f=p_Xg$.  Hence $\C(P,P_i)\ne0$ for some $i$, and
\[
 \{P\in\ind\Pcal\mid\C(P,X)\ne0\}
 \subseteq
 \bigcup_{i=1}^t\{P\in\ind\Pcal\mid\C(P,P_i)\ne0\}.
\]
The right side is finite because $\Pcal$ is locally finite.

Choose $Q_X$ as in (D$_-$), and an $\E$-triangle
\[
 \Omega Q_X\longrightarrow P_X\longrightarrow Q_X\dashrightarrow,
~~\mbox{where}~P_X\in\Pcal.
\]
Set
\[
 r_X^-=[\Omega Q_X]+[Q_X]-[P_X]\in\Rel(\C).
\]
Since $\Rel(\C)=\AR(\C)$, write
\[
 r_X^-=\sum_{j=1}^m a_j[\delta_j]
\]
with $a_j\in\mathbb Z$ and $\delta_j$ Auslander--Reiten $\E$-triangles.
Then
\[
 \Supp\chi^-(r_X^-)
 \subseteq\bigcup_{j=1}^m\Supp\chi^-([\delta_j]),
\]
which is finite by (EFS). By Lemma~\ref{lem:stable-defect}, it yields
\[
 \{U\in\ind\C\mid\underline\C(U,Q_X)\ne0\}
 \subseteq\Supp\chi^-(r_X^-).
\]
Using (D$_-$), we obtain
\[
 \Supp\C(-,X)
 \subseteq
 \{P\in\ind\Pcal\mid\C(P,X)\ne0\}
 \cup
 \{U\in\ind\C\mid\underline\C(U,Q_X)\ne0\},
\]
so $\Supp\C(-,X)$ is finite.

To prove finiteness of the outgoing support, choose an $\E$-triangle
\[
 X\xrightarrow{i_X}I^X\longrightarrow\Sigma X\dashrightarrow,
 ~~\mbox{where}~I^X\in\Ical,
\]
and write $I^X=I_1\oplus\cdots\oplus I_s$ with $I_j\in\ind\Ical$.  If
$I\in\ind\Ical$ and $0\ne f\colon X\to I$, injectivity of $I$ gives a
morphism $h\colon I^X\to I$ such that $f=hi_X$.  Hence
$\C(I_j,I)\ne0$ for some $j$, and therefore
\[
 \{I\in\ind\Ical\mid\C(X,I)\ne0\}
 \subseteq
 \bigcup_{j=1}^s
 \{I\in\ind\Ical\mid\C(I_j,I)\ne0\}.
\]
The right side is finite because $\Ical$ is locally finite.

Choose $T_X$ as in (D$_+$), and choose an $\E$-triangle
\[
 T_X\longrightarrow I_X\longrightarrow\Sigma T_X\dashrightarrow,
 ~~\mbox{where}~ I_X\in\Ical.
\]
Set
\[
 r_X^+=[T_X]+[\Sigma T_X]-[I_X]\in\Rel(\C).
\]
Since $\Rel(\C)=\AR(\C)$, write
\[
 r_X^+=\sum_{j=1}^n b_j[\varepsilon_j]
\]
with $b_j\in\mathbb Z$ and $\varepsilon_j$ Auslander--Reiten
$\E$-triangles.  By (EFS),
\[
 \Supp\chi^+(r_X^+)
 \subseteq\bigcup_{j=1}^n\Supp\chi^+([\varepsilon_j])
\]
is finite.  The second formula in Lemma~\ref{lem:stable-defect} gives
\[
 \{U\in\ind\C\mid\overline\C(T_X,U)\ne0\}
 \subseteq\Supp\chi^+(r_X^+).
\]
Using (D$_+$), we obtain
\[
 \Supp\C(X,-)
 \subseteq
 \{I\in\ind\Ical\mid\C(X,I)\ne0\}
 \cup
 \{U\in\ind\C\mid\overline\C(T_X,U)\ne0\},
\]
so $\Supp\C(X,-)$ is finite.  Thus $\C$ is locally finite.
\end{proof}

\begin{remark}\label{rem:role-hypotheses}
The proof separates the two sides of local finiteness.  To obtain
$(\mathrm{LF}_-)$, it uses enough projectives, the local finiteness of
$\Pcal$, the finiteness of $\Supp\chi^-([\rho])$ for
Auslander--Reiten $\E$-triangles $\rho$, the condition (D$_-$), and
$\Rel(\C)=\AR(\C)$.  The dual datas give $(\mathrm{LF}_+)$.  Thus the
projective and injective boundary assumptions, the Euler functions, and the
detection conditions enter at different stages of the argument.

The boundary assumption cannot simply be omitted.  In
Example~\ref{ex:one-sided-not-self-dual}, the extriangulated structure is
split, so every object is projective and injective.  Hence the category has
enough projectives and injectives, (EFS) is vacuous, and
(D$_-$) and {(D$_+$)} hold automatically.  Moreover, every
$\E$-triangle is split, so
\[
 \Rel(\C)=\AR(\C)=0.
\]
Nevertheless, the category is not locally finite.  In this example
$\Pcal=\Ical=\C$, and these boundary subcategories are not locally
finite.  Thus the relation equality, even together with the Euler and
detection conditions, does not by itself control morphisms through
projective or injective objects.

We do not claim that the remaining assumptions are logically independent in full generality.  The criteria below are useful because these assumptions arise naturally in many familiar settings.  Exact structures and Serre duality give (EFS).  Weak cogenerators and weak generators give the detection conditions.  Projective generators and injective cogenerators ensure that the corresponding boundary subcategories are of finite type.
\end{remark}

By Lemma~\ref{lem:weak-detectors}, the weak detector conditions provide a useful way to verify
{(D$_-$)} and {(D$_+$)}. We first record a direct consequence of the main theorem before giving more concrete criteria.

\begin{corollary}\label{cor:one-sided-equivalence}
Let $\C$ be a skeletally small Hom-finite Krull--Schmidt $R$-linear
extriangulated category with enough projectives and enough injectives.  Assume
that $\Pcal$ and $\Ical$ are locally finite and that $\C$ satisfies
\emph{(EFS)}, \emph{(D$_-$)} and \emph{(D$_+$)}.  Then the following are
equivalent:
\begin{enumerate}[label=\textup{(\arabic*)},leftmargin=2.7em]
\item $(\mathrm{LF}_-)$ holds;
\item $(\mathrm{LF}_+)$ holds;
\item $\C$ is locally finite;
\item $\Rel(\C)=\AR(\C)$.
\end{enumerate}
\end{corollary}

\begin{proof}
Either \textup{(1)} or \textup{(2)} implies \textup{(4)} by
Theorem~\ref{thm:one-sided-relations}.  Under the standing structural
assumptions, \textup{(4)} implies \textup{(3)} by
Theorem~\ref{thm:main}.  Finally, \textup{(3)} implies both
\textup{(1)} and \textup{(2)} by definition.
\end{proof}

\begin{corollary}[A generator and cogenerator criterion]\label{cor:generator-cogenerator}
Let $\C$ be a skeletally small Hom-finite Krull--Schmidt $R$-linear
extriangulated category which has enough projectives add$P$ and enough injectives add$I$.  Assume that $\C$ satisfies
\emph{(EFS)}, that $\underline\C$ has a weak cogenerator, and that
$\overline\C$ has a weak generator.  Then
\[
 \C\text{ is of finite type }
 ~~\Longleftrightarrow~~
 \Rel(\C)=\AR(\C).
\]
\end{corollary}

\begin{proof}
By assumption, $\C$ has enough projectives given by add$P$ and enough injectives given by add$I$.
By Lemma~\ref{lem:weak-detectors}, {(D$_-$)} and {(D$_+$)} hold.
Since $\C$ has projective generator $P$, the implication from right to left follows from
Theorem~\ref{thm:main}.  The converse follows from
Theorem~\ref{thm:one-sided-relations}.
\end{proof}

We apply it to $0$-Auslander extriangulated categories in the sense of
{\rm\cite[Definition~3.7]{GorskyNakaokaPalu2023}}.

\begin{corollary}\label{cor:zero-Auslander}
Let $\C$ be a skeletally small Hom-finite Krull--Schmidt $R$-linear
$0$-Auslander extriangulated category. Assume that $\Pcal=\text{add} P$ and that $\C$ satisfies \emph{(EFS)}.  Then
\[
 \C\text{ is of finite type}
 ~~\Longleftrightarrow~~
 \Rel(\C)=\AR(\C).
\]
\end{corollary}

\begin{proof}
A $0$-Auslander extriangulated category has enough projectives and enough
injectives.  Since $\Pcal$ is of finite type, $\Ical$ is also of finite type, they are locally
finite.  By \cite[Proposition~3.11]{GorskyNakaokaPalu2023}, the costable
category $\overline\C=\C/[\Ical]$ is equivalent to
$\modcat\overline{\Pcal}$, where $\overline{\Pcal}$ is the image of
$\Pcal$ modulo the projective-injective objects.  Since $\Pcal$ is of finite
type, $\overline{\Pcal}$ has an additive generator $P_0$.  Let $G$ be the
object of $\overline\C$ corresponding to the representable functor
$\overline{\Pcal}(-,P_0)$.  We claim that $G$ is a weak generator of
$\overline\C$.  Indeed, suppose that $X\in\C$ satisfies
$\overline\C(G,X)=0$.  By the above equivalence and Yoneda, we have
$$
 0=\Hom_{\modcat\overline{\Pcal}}
   (\overline{\Pcal}(-,P_0),F_X)
 \simeq F_X(P_0),
$$
where $F_X$ denotes the functor corresponding to $X$.
 Since every object of
$\overline{\Pcal}$ belongs to $\add P_0$, it follows that $F_X=0$.
Hence $X$ is zero in $\overline\C$, so $1_X$ factors through an injective
object. Thus $X$ is injective, and the claim follows.

Now apply the same proposition to $\C^{\mathrm{op}}$, which is again
$0$-Auslander.  Its costable category is $(\underline\C)^{\mathrm{op}}$.
The projective objects of $\C^{\mathrm{op}}$ are precisely the injective
objects of $\C$.  Hence the finite type assumption on $\Ical$ gives a weak
generator of $(\underline\C)^{\mathrm{op}}$ by the same argument.  Viewed in
$\underline\C$, this object is a weak cogenerator.  Indeed, if
$\underline\C(X,Q)=0$, then $X$ is zero in $\underline\C$ and therefore
projective.

By Lemma~\ref{lem:weak-detectors}, we obtain {(D$_-$)} and
{(D$_+$)}.  If $\Rel(\C)=\AR(\C)$, local finiteness follows from
Theorem~\ref{thm:main}. Hence $\C$ is of finite type. The converse follows from
Theorem~\ref{thm:one-sided-relations}.
\end{proof}

We next give an example in which the hypotheses of
Corollary~\ref{cor:zero-Auslander} arise naturally from a relative
extriangulated structure that is generally neither exact nor triangulated.

\begin{example}\label{ex:relative-cluster}
Let $k$ be a field and let $\mathcal T$ be a skeletally small Hom-finite
Krull--Schmidt $k$-linear $2$-Calabi--Yau triangulated category.  Let
$T$ be a cluster tilting object.  We endow the underlying additive category of
$\mathcal T$ with the relative extriangulated structure whose $\E$-triangles
are the triangles
$$
 X\longrightarrow Y\longrightarrow Z\xrightarrow{~\eta~}\Sigma X
$$
for which $\eta$ factors through an object of $\add\Sigma T$.  We denote the
resulting extriangulated category by $\mathcal T_T$.

By \cite[Section~3.3.6]{GorskyNakaokaPalu2023}, $\mathcal T_T$ is a
$0$-Auslander extriangulated category with
\[
 \Proj\mathcal T_T=\add T,
~~
 \Inj\mathcal T_T=\add\Sigma^2 T.
\]
Hence its projective and injective subcategories are of finite type.  Since
$\mathcal T$ is $2$-Calabi--Yau, $\Sigma^2$ is a Serre functor on the
underlying additive category, and Proposition~\ref{prop:Serre-EFS} shows
that $\mathcal T_T$ satisfies {(EFS)}.
By Corollary~\ref{cor:zero-Auslander}, we give
\[
 \mathcal T\text{ is of finite type}
 ~~\Longleftrightarrow~~
 \Rel(\mathcal T_T)=\AR(\mathcal T_T).
\]

For a Dynkin cluster category, the condition on the left is automatic.  Thus
this includes the familiar relative examples of finite type.  To compare our
result with that of Zhu and Zhuang \cite{ZhuZhuang2021}, note that
$K_0(\mathcal T_T)$ is precisely their relative Grothieck group.  Its defining
relations come from those triangles whose connecting morphism factors through
$\add\Sigma T$.  Moreover, Zhu and Zhuang
\cite[Lemma~4.13 and the proof of Theorem~4.14]{ZhuZhuang2021} showed that the
Auslander--Reiten $\E$-triangles for this relative structure are exactly the
ambient Auslander--Reiten triangles satisfying the same factorization
condition.  Hence the equivalence above gives another proof of the cluster
tilting object and $2$-Calabi--Yau special case of their relative criterion
\cite[Theorem~4.14]{ZhuZhuang2021}.

\end{example}

\section{Applications and examplesl cases}\label{sec:applications}

\hspace{-0.5mm}{\bf {\thesection.1} Module categories over Artin algebras.}~
Let $\Lambda$ be an Artin $R$-algebra and let $\modcat\Lambda$ be the
category of finitely generated left $\Lambda$-modules.  Choose
representatives $S_1,\dots,S_t$ of the simple left modules and put
\[
 S_\Lambda=S_1\oplus\cdots\oplus S_t.
\]

\begin{lemma}\label{lem:module-detectors}
The object $S_\Lambda$ is a weak cogenerator of
$\underline{\modcat}\Lambda$ and a weak generator of
$\overline{\modcat}\Lambda$.  Equivalently,
\begin{align*}
 {\rm\underline{\emph{\Hom}}}_\Lambda(U,S_\Lambda)=0
 &~~\Longrightarrow~~U\text{ is projective},\\[2mm]
 {\rm\overline{ \emph{\Hom}}}_\Lambda(S_\Lambda,U)=0
 &~~\Longrightarrow~~ U\text{ is injective}.
\end{align*}
\end{lemma}

\proof The first assertion is \cite[Proposition~4.7]{Enomoto2019}, since
$\Lambda/\rad\Lambda$ and $S_\Lambda$ have the same indecomposable simple
summands.  The second follows from the same result for
$\Lambda^{\mathrm{op}}$ and the standard Artin duality.  \qed

\vspace{2mm}
Butler \cite{Butler1981} proved that finite representation type implies that the relations of
the Grothendieck group are generated by almost split sequences.  Auslander
\cite{Auslander1984} proved the converse.  Since $\modcat\Lambda$ is locally
finite if and only if $\Lambda$ is of finite representation type, we obtain
the following formulation.

\begin{theorem}\label{thm:module}{\rm\cite{Butler1981,Auslander1984}}
Let $R$ be a commutative artin ring and let $\Lambda$ be an Artin
$R$-algebra.  The following are equivalent.
\begin{enumerate}[label=\textup{(\arabic*)}]
\item $\Lambda$ is of finite representation type;
\item $\modcat\Lambda$ is locally finite;
\item $\Rel(\modcat\Lambda)=\AR(\modcat\Lambda)$.
\end{enumerate}
\end{theorem}


More general exact categories related to the Auslander--Butler theorem are
studied in \cite{Enomoto2019}.
\vspace{1mm}

\hspace{-5mm}{\bf {\thesection.2} Triangulated categories.}~
Let $\C$ be a Hom-finite Krull--Schmidt triangulated category with suspension
$[1]$, regarded as an extriangulated category via
$\E(X,Y)=\C(X,Y[1])$.  Its only projective and injective object is zero, and
it has enough projectives and injectives in the extriangulated sense; for
example, one may use the trivial triangles
\[
 X[-1]\longrightarrow0\longrightarrow X\xrightarrow{1_X}X~~\mbox{and}~~
 X\longrightarrow0\longrightarrow X[1]\xrightarrow{1_{X[1]}}X[1].
\]

\begin{corollary}\label{cor:triangulated}
Let $\C$ be a skeletally small Hom-finite Krull--Schmidt $R$-linear
triangulated category.  The following conditions are equivalent:
\begin{enumerate}[label=\textup{(\arabic*)},leftmargin=2.7em]
\item $|\Supp\C(-,X)|<\infty$ for every $X\in\ind\C$;
\item $|\Supp\C(X,-)|<\infty$ for every $X\in\ind\C$;
\item $\C$ is locally finite;
\item $\Rel(\C)=\AR(\C)$.
\end{enumerate}
\end{corollary}

\begin{proof}
Here $\Pcal=\Ical=0$ and $\underline\C=\overline\C=\C$.  Thus (D$_-$)
and (D$_+$) are automatic, while (EFS) follows from
\eqref{eq:tri-support}. The result now follows from
Corollary~\ref{cor:one-sided-equivalence}.
\end{proof}

For $R=k$ a field, the implication \textup{(3)}$\Rightarrow$\textup{(4)}
recovers \cite[Theorem~2.1]{XiaoZhu2002}, while the converse
\textup{(4)}$\Rightarrow$\textup{(3)} recovers
\cite[Theorem~3.14 and Remark~3.15]{Zhou2021}.  The equivalence of
\textup{(1)} and \textup{(2)} shows that one-sided local finiteness is
self-dual in the triangulated setting.

\vspace{3mm}

\hspace{-5mm}{\bf {\thesection.3} Frobenius extriangulated categories.}~\label{subsec:frobenius}
Let $\C$ be a Frobenius extriangulated category and write
\[
 \underline{\C}=\C/[\Pcal],
 ~~\Pcal=\Proj\C=\Inj\C.
\]
By \cite[Corollary~7.4 and Remark~7.5]{NakaokaPalu2019}, the stable category
$\underline{\C}$ is triangulated.  It is also Krull--Schmidt.  Indeed, after removing
projective summands, every object of $\C$ remains a finite direct sum of
indecomposable objects.  The endomorphism ring of the image of a
nonprojective indecomposable is a quotient of a local ring and is therefore
local.  For Frobenius exact categories, the corresponding statement on
relation subgroups appears in Haugland
\cite[Lemma~3.2]{Haugland2022}.  We now record the extriangulated version.

\begin{proposition}\label{prop:frobenius-stable-relations}
Let $\C$ be a skeletally small Hom-finite Krull--Schmidt $R$-linear
Frobenius extriangulated category.  If
\[
 \Rel(\C)=\AR(\C),
\]
then
\[
 \Rel(\underline{\C})=\AR(\underline{\C}),
\]
and the triangulated category $\underline{\C}$ is locally finite.
\end{proposition}

\begin{proof}
Let $\pi\colon\C\to\underline{\C}$ be the quotient functor.  It induces a
surjective homomorphism
\[
 \pi_*\colon K_0^{\mathrm{sp}}(\C)
 \longrightarrow K_0^{\mathrm{sp}}(\underline{\C}),
 ~~[X]\longmapsto[{X}].
\]
By the standard Frobenius triangulation on $\underline{\C}$
\cite[Corollary~7.4 and Remark~7.5]{NakaokaPalu2019}, the quotient functor
sends every $\E$-triangle in $\C$ to a distinguished triangle in
$\underline{\C}$, and, by construction, every distinguished triangle is
isomorphic to one obtained in this way.  For the explicit construction under
(WIC), compare \cite[Definition~6.15 and Theorem~6.20]{NakaokaPalu2019}.
Hence
\[
 \pi_*(\Rel(\C))=\Rel(\underline{\C}).
\]
By \cite[Proposition~5.11]{IyamaNakaokaPalu2024}, an almost split extension
remains almost split in the ideal quotient by any full additive subcategory
of projective-injective objects.  Taking this subcategory to be $\Pcal$, the
image of an Auslander--Reiten $\E$-triangle in $\C$ is an
Auslander--Reiten triangle in $\underline{\C}$, and therefore
\[
 \pi_*(\AR(\C))\subseteq\AR(\underline{\C}).
\]
The hypothesis gives
\[
 \Rel(\underline{\C})
 =\pi_*(\Rel(\C))
 =\pi_*(\AR(\C))
 \subseteq\AR(\underline{\C})
 \subseteq\Rel(\underline{\C}),
\]
so $\Rel(\underline{\C})=\AR(\underline{\C})$.
The stable category is skeletally small and Hom-finite, and it is
Krull--Schmidt as noted above. By Corollary~\ref{cor:triangulated}, we  show
that $\underline{\C}$ is locally finite.
\end{proof}

Notice that no finiteness assumption on the indecomposable projective
objects of $\C$ is needed here.  The conclusion concerns the stable category
$\underline{\C}$, where the projective and injective boundary has
vanished.


The examples below illustrate several features of Theorem~\ref{thm:main} that
do not appear in the classical cases.  We verify local finiteness and the
existence of enough projectives and injectives directly.  The remaining
conditions then follow from the implication
\textup{(1)}$\Rightarrow$\textup{(2)} of Theorem~\ref{thm:main}.  These
examples include categories of infinite type, genuinely extriangulated
structures, and cases in which the projective and injective subcategories
contain infinitely many indecomposable objects.

\begin{example}[Two-term projective complexes]\label{ex:two-term}
Let $k$ be a field and let $Q$ be a Dynkin quiver.  Put $\Lambda=kQ$ and consider
$$
 \C=K^{[-1,0]}(\proj\Lambda)
 \subseteq K^b(\proj\Lambda)
$$
with the extriangulated structure induced from the ambient triangulated
category.  Categories of complexes of fixed size were studied in
\cite{BautistaSoutoSalorioZuazua2005}.  In the present homotopy--extriangulated
setting, it follows from \cite[Section~3.3.3]{GorskyNakaokaPalu2023} that this
is a reduced $0$-Auslander extriangulated category.  In particular, it has
enough projectives and injectives.  Its projective objects are the stalk
complexes concentrated in degree $0$, while its injective objects are the
stalk complexes concentrated in degree $-1$.

Since $\Lambda$ is hereditary, every object $X$ of $\C$, viewed in
$\Db(\modcat\Lambda)$, is isomorphic to
$$
 H^{-1}(X)[1]\oplus H^0(X),
$$
and $H^{-1}(X)$ is projective.  See, for example,
\cite[Chapter~I, Section~5]{Happel1988}.  Under the fully faithful equivalence
$K^b(\proj\Lambda)\simeq\Db(\modcat\Lambda)$, the module $H^0(X)$ is
represented by its projective resolution.  This resolution is concentrated in degrees
$-1$ and $0$.  Thus the indecomposable objects of $\C$ are represented by
indecomposable projective modules in degree $-1$ or by projective resolutions
of indecomposable modules.  Since $Q$ is Dynkin, there are only finitely many
such objects.  Hence $\C$ is of finite type and therefore locally finite.

The induced extriangulated structure is neither exact nor triangulated.
For a nonzero projective module $P$, the triangle
$$
 P\longrightarrow0\longrightarrow P[1]\longrightarrow P[1]
$$
has all three terms in $K^{[-1,0]}(\proj\Lambda)$.  Hence $P\to0$ is an
inflation in the induced extriangulated structure.  However, it is not a
monomorphism.  Thus the structure is not exact.  It is not triangulated
either because it has nonzero projective objects.
\end{example}

The preceding example is of finite type.  We next give an example of
infinite type which is neither exact nor triangulated.

Let $k$ be a field and let
\[
 Q_n:\quad 1\longrightarrow2\longrightarrow\cdots\longrightarrow n,
 ~ n\ge3,
\]
and set $\Lambda_n=kQ_n$.  Let $\mathcal B_n$ be the full subcategory of
$\modcat\Lambda_n$ consisting of modules whose component at the vertex $1$
is zero.  Thus $\mathcal B_n$ is a Serre subcategory, equivalent to the
module category of the linearly oriented quiver on the vertices
$2,\cdots,n$.  Define
\[
 \mathcal X_n=
 \{X\in\Db(\modcat\Lambda_n)\mid
   {\rm H}^r(X)\in\mathcal B_n\text{ for every }r\ne0\}.
\]

\begin{proposition}\label{prop:infinite-genuine-example}
The category $\mathcal X_n$ is an extension-closed Hom-finite
Krull--Schmidt subcategory of $\Db(\modcat\Lambda_n)$, hence it is
extriangulated.  It has the following properties.
\begin{enumerate}[label=\textup{(\roman*)},leftmargin=2.7em]
\item $\mathcal X_n$ is locally finite but not of finite type;
\item $\mathcal X_n$ has enough projectives and injectives, and
\[
 \Proj\mathcal X_n=\add P_1,
~~
 \Inj\mathcal X_n=\add S_1,
\]
where $P_1$ is the indecomposable projective at the vertex $1$ and $S_1$ is
the simple module at the vertex $1$;
\item $\mathcal X_n$ is neither exact nor triangulated.
\end{enumerate}
\end{proposition}

\begin{proof}
Since $\mathcal B_n$ is a Serre subcategory, the long exact cohomology
sequence associated with a triangle shows that $\mathcal X_n$ is extension
closed.  Indeed, if $X,Z\in\mathcal X_n$ and
$X\to Y\to Z\to X[1]$ is a triangle, then for $r\ne0$ there is a short
exact sequence
\[
 0\longrightarrow
 \operatorname{coker}\bigl(H^{r-1}(Z)\to H^r(X)\bigr)
 \longrightarrow H^r(Y)\longrightarrow
 \ker\bigl(H^r(Z)\to H^{r+1}(X)\bigr)
 \longrightarrow0.
\]
Both end terms belong to $\mathcal B_n$, and hence so does $H^r(Y)$.  Thus
$\mathcal X_n$ inherits an extriangulated structure from
$\Db(\modcat\Lambda_n)$ by \cite[Remark~2.18]{NakaokaPalu2019}.

Write $M_{i,j}$ for the interval module supported on the vertices
$i,i+1,\cdots,j$; these are precisely the indecomposable
$\Lambda_n$-modules, see \cite{ARS1995}.  Since $\Lambda_n$ is hereditary,
every indecomposable object of $\Db(\modcat\Lambda_n)$ is a shift of an
indecomposable module; see \cite[Chapter~I, Section~5]{Happel1988}.  Hence
\[
 \ind\mathcal X_n=
 \{M_{i,j}\mid1\le i\le j\le n\}
 \cup
 \{M_{i,j}[r]\mid2\le i\le j\le n,\ r\ne0\}.
\]
In particular, the objects $S_2[r]$, $r\in\mathbb Z$, are pairwise
nonisomorphic, so $\mathcal X_n$ is not of finite type.  Local finiteness can
be checked directly.  For indecomposable modules $M,N$ and integers $r,s$,
\[
 \Hom_{\Db(\modcat\Lambda_n)}(M[r],N[s])
 \simeq \Ext^{s-r}_{\Lambda_n}(M,N),
\]
and the right-hand side can be nonzero only for $s-r=0$ or $1$, since
$\Lambda_n$ is hereditary.  There are only finitely many interval modules
$M_{i,j}$.  Hence, for any fixed indecomposable object of $\mathcal X_n$,
only finitely many indecomposable objects of $\mathcal X_n$ can occur in
either Hom support.  Thus $\mathcal X_n$ is locally finite.

We determine the projective objects.  Put $P_1=M_{1,n}$.  If
$N\in\mathcal B_n$, then
\[
 \Hom_{\Lambda_n}(P_1,N)=0.
\]
Let $Y$ be an indecomposable object of $\mathcal X_n$.  If $Y$ is a module
concentrated in degree $0$, then
\[
 \Hom_{\Db(\modcat\Lambda_n)}(P_1,Y[1])
 =\Ext^1_{\Lambda_n}(P_1,Y)=0.
\]
If $Y=N[r]$ with $N\in\mathcal B_n$ and $r\ne0$, then projectivity of
$P_1$ and the hereditary property show that
$\Hom(P_1,N[r+1])$ can be nonzero only when $r=-1$, in which case it is
$\Hom_{\Lambda_n}(P_1,N)=0$.  Therefore
\[
 \Hom_{\Db(\modcat\Lambda_n)}(P_1,Y[1])=0
 ~~(Y\in\mathcal X_n),
\]
so $P_1$ is projective in $\mathcal X_n$.  Every nonzero object
$N[r]$ with $N\in\mathcal B_n$ is not projective, since
\[
 \E_{\mathcal X_n}(N[r],N[r-1])
 =\End_{\Db(\Lambda_n)}(N[r])\ne0.
\]
For $1\le j<n$, the nonsplit exact sequence
\[
 0\longrightarrow M_{j+1,n}\longrightarrow P_1
 \longrightarrow M_{1,j}\longrightarrow0
\]
shows that $M_{1,j}$ is not projective.  Hence
$\Proj\mathcal X_n=\add P_1$.

Dually, put $S_1=M_{1,1}$.  The module $S_1$ is injective over
$\Lambda_n$, and $\Hom_{\Lambda_n}(N,S_1)=0$ for
$N\in\mathcal B_n$.  If $Y$ is a module concentrated in degree $0$, then
$\Hom(Y,S_1[1])=\Ext^1_{\Lambda_n}(Y,S_1)=0$.  If
$Y=N[r]$ with $N\in\mathcal B_n$ and $r\ne0$, injectivity of $S_1$ and
the hereditary property show that $\Hom(N[r],S_1[1])$ can be nonzero only
when $r=1$, in which case it is $\Hom_{\Lambda_n}(N,S_1)=0$.
Thus $S_1$ is injective in $\mathcal X_n$.  Every nonzero
$N[r]$ with $N\in\mathcal B_n$ is not injective because
\[
 \E_{\mathcal X_n}(N[r+1],N[r])
 =\End_{\Db(\Lambda_n)}(N[r+1])\ne0.
\]
If $j>1$, then
\[
 \E_{\mathcal X_n}(S_j[1],M_{1,j})
 =\Hom_{\Lambda_n}(S_j,M_{1,j})\ne0,
\]
so $M_{1,j}$ is not injective.  Thus
$\Inj\mathcal X_n=\add S_1$.

We finally check that there are enough projectives and injectives.  For
$N\in\mathcal B_n$ and $r\in\mathbb Z$, the trivial triangle
\[
 N[r-1]\longrightarrow0\longrightarrow N[r]
 \xrightarrow{1}N[r]
\]
provides a projective approximation.  For $M_{1,j}$ one uses the triangle
induced by
\[
 0\longrightarrow M_{j+1,n}\longrightarrow P_1
 \longrightarrow M_{1,j}\longrightarrow0,
\]
with the convention $M_{n+1,n}=0$.  Thus $\mathcal X_n$ has enough
projectives.  Similarly, the trivial triangle gives an injective
approximation for $N[r]$, while the exact sequence
\[
 0\longrightarrow M_{2,j}\longrightarrow M_{1,j}
 \longrightarrow S_1\longrightarrow0
\]
yields, after rotation, a triangle
\[
 M_{1,j}\longrightarrow S_1\longrightarrow M_{2,j}[1]
 \longrightarrow M_{1,j}[1],
\]
with the convention $M_{2,1}=0$.  Hence $\mathcal X_n$ has enough
injectives.

The category is not triangulated because it has the nonzero projective object
$P_1$.  It is not exact either.  For example, consider the rotated triangle
$$
 M_{1,2}\longrightarrow S_1\longrightarrow S_2[1]
 \longrightarrow M_{1,2}[1]
$$
in $\mathcal X_n$.  The morphism $M_{1,2}\to S_1$ is an inflation, but it is not a monomorphism.
\end{proof}

\begin{example}\label{ex:graded-dual-numbers}
Let $k$ be a field and let
$
 \Lambda=k[\varepsilon]/(\varepsilon^2)
$
be the algebra of dual numbers, regarded as a $\mathbb Z$-graded algebra
with $\varepsilon$ homogeneous of degree $1$.  Let
$\C=\operatorname{grmod}\Lambda$ be the category of finite-dimensional
$\mathbb Z$-graded left $\Lambda$-modules with morphisms of degree zero.
Then $\C$ is a Hom-finite Krull--Schmidt exact category which is locally
finite but not of finite type.  Moreover,
\[
 |\ind\Proj\C|=|\ind\Inj\C|=\infty,
\]
and $\C$ satisfies all the hypotheses of Theorem~\ref{thm:main}(2).
\end{example}

\begin{proof}
Put
$
 S=\Lambda/(\varepsilon).
$
For a graded module
\[
 M=\bigoplus_{i\in\mathbb Z}M_i,
\]
multiplication by $\varepsilon$ gives $k$-linear maps
$
 M_i\xrightarrow{\ \varepsilon\ }M_{i+1}.
$
Since $\varepsilon^2=0$, the composite of any two consecutive maps is
zero.  Thus $M$, together with multiplication by $\varepsilon$, determines
a bounded complex of finite-dimensional $k$-vector spaces
\[
 \cdots\longrightarrow M_{i-1}
 \xrightarrow{\ \varepsilon\ }M_i
 \xrightarrow{\ \varepsilon\ }M_{i+1}
 \longrightarrow\cdots.
\]
Conversely, every such bounded complex defines a finite-dimensional graded
$\Lambda$-module by letting $\varepsilon$ act as the differential.
Moreover, morphisms of degree zero correspond precisely to morphisms of
complexes.

Every bounded complex of finite-dimensional vector spaces is isomorphic
to a direct sum of stalk complexes and two-term contractible complexes.
Indeed, for each $i$ one may choose a decomposition
\[
 M_i=B_i\oplus H_i\oplus C_i
\]
such that
\[
 \ker(\varepsilon\colon M_i\to M_{i+1})=B_i\oplus H_i,
~~
 \varepsilon\colon C_i\xrightarrow{\sim}B_{i+1}.
\]
It follows that the indecomposable objects of $\C$ are precisely
\[
 S\langle n\rangle
 ~~\text{and}~~
 \Lambda\langle n\rangle,
 ~~ n\in\mathbb Z,
\]
where $\langle n\rangle$ denotes the grading shift.  In particular,
the objects $S\langle n\rangle$, $n\in\mathbb Z$, are pairwise
nonisomorphic, and hence $\C$ is not of finite type.

We next prove local finiteness.  The module $S\langle n\rangle$ is
concentrated in one degree, whereas $\Lambda\langle n\rangle$ is
concentrated in two consecutive degrees.  A morphism in $\C$ has degree
zero.  Hence, for a fixed indecomposable object $X$, a nonzero morphism
$
 Y\longrightarrow X
 ~~\text{or}~~
 X\longrightarrow Y
$
is possible only if the graded supports of $X$ and $Y$ intersect.
Since there are only two types of indecomposable objects for each grading
shift, this leaves only finitely many possibilities for $Y$.  Therefore
both
\[
 \Supp\C(-,X)
~~\text{and}~~
 \Supp\C(X,-)
\]
are finite for every $X\in\ind\C$.  Thus $\C$ is locally finite.

We now determine the projective and injective objects.  For each
$n\in\mathbb Z$, set
\[
 S_n=S\langle n\rangle,
 ~~
 P_n=\Lambda\langle n\rangle.
\]
Since $P_n$ is a grading shift of the free module $\Lambda$, it is
projective.  Moreover, there is a nonsplit short exact sequence
\[
 0\longrightarrow S_{n+1}
 \longrightarrow P_n
 \longrightarrow S_n
 \longrightarrow0.
\]
In particular, $S_n$ is not projective.  Hence
$
 \ind\Proj\C=\{P_n\mid n\in\mathbb Z\}.
$
On the other hand, $\Lambda$ is a graded Frobenius algebra.  Indeed, the
linear functional
\[
 \lambda\colon\Lambda\longrightarrow k,
~~
 \lambda(1)=0,~~
 \lambda(\varepsilon)=1,
\]
defines an associative bilinear form
$
 \Lambda\times\Lambda\longrightarrow k,
~~
 (a,b)\longmapsto\lambda(ab),
$
which is nondegenerate.  Thus $D\Lambda$ is isomorphic, as a graded
$\Lambda$-module, to a grading shift of $\Lambda$.
Consequently the
grading shifts $P_n$ are injective as well as projective.  Since
\[
 0\longrightarrow S_n
 \longrightarrow P_{n-1}
 \longrightarrow S_{n-1}
 \longrightarrow0
\]
is nonsplit, $S_n$ is not injective.  Therefore
$
 \ind\Inj\C=\{P_n\mid n\in\mathbb Z\}.
$
In particular,
\[
 |\ind\Proj\C|
 =
 |\ind\Inj\C|
 =
 \infty.
\]
The preceding short exact sequences also show directly that $\C$ has
enough projectives and enough injectives.

For completeness, we verify the equality of relation subgroups directly.
For every $n\in\mathbb Z$, the short exact sequence
\[
 0\longrightarrow S_{n+1}
 \longrightarrow P_n
 \longrightarrow S_n
 \longrightarrow0
\]
is an Auslander--Reiten sequence.  To see this, note first that the
epimorphism $P_n\to S_n$ is right almost split.  Indeed, using the above
classification of indecomposables, a nonzero morphism from an
indecomposable object to $S_n$ which is not an isomorphism can only come
from $P_n$, and hence factors through $P_n\to S_n$.  Dually,
$S_{n+1}\to P_n$ is left almost split.

Since $\C$ is a length category and its simple objects are precisely the
$S_n$, one has
\[
 K_0(\C)\simeq
 \bigoplus_{n\in\mathbb Z}\mathbb Z[S_n].
\]
On the other hand,
\[
 K_0^{\mathrm{sp}}(\C)
 =
 \bigoplus_{n\in\mathbb Z}
 \bigl(\mathbb Z[S_n]\oplus\mathbb Z[P_n]\bigr).
\]
The above short exact sequence gives
\[
 [P_n]=[S_n]+[S_{n+1}]
 ~~\text{in }K_0(\C).
\]
Therefore the kernel of the canonical homomorphism
\[
 K_0^{\mathrm{sp}}(\C)\longrightarrow K_0(\C)
\]
is generated by the elements
\[
 [S_n]+[S_{n+1}]-[P_n],
~~n\in\mathbb Z.
\]
These are exactly the relations of the Auslander--Reiten sequences above.
Consequently,
\[
 \Rel(\C)=\AR(\C).
\]
Finally, $\C$ is abelian and hence exact, and it is Hom-finite and
Krull--Schmidt.  Since $\C$ is locally finite and has enough projectives
and injectives, Lemma~\ref{lem:locally-finite-detection} shows that
$\Proj\C$ and $\Inj\C$ are locally finite and that $\C$ satisfies
{(EFS)}, {(D$_-$)} and {(D$_+$)}.  Thus every hypothesis
of Theorem~\ref{thm:main} is satisfied.
\end{proof}


\section{A criterion using Serre duality}\label{sec:serre}

Ordinary Serre duality interchanges the two Hom support conditions.  Thus,
once one side is known to be finite, no separate argument is needed for the
opposite projective or injective boundary.

\begin{theorem}\label{thm:Serre-one-sided}
Let $\C$ be a skeletally small Hom-finite Krull--Schmidt $R$-linear
extriangulated category admitting a Serre functor $S$, and assume that
\[
 \Rel(\C)=\AR(\C).
\]
Suppose that one of the following two sets of conditions holds:
\begin{enumerate}[label=\textup{(\arabic*)},leftmargin=2.7em]
\item $\C$ has enough projectives, $\Pcal$ is locally finite, and $\C$
satisfies \emph{(D$_-$)};
\item $\C$ has enough injectives, $\Ical$ is locally finite, and $\C$
satisfies \emph{(D$_+$)}.
\end{enumerate}
Then $\C$ is locally finite.
\end{theorem}

\begin{proof}
We prove the assertion under \textup{(1)}, the proof under \textup{(2)} is
dual.  By Proposition~\ref{prop:Serre-EFS}, we show that
$\Supp\chi^-([\delta])$ is finite for every Auslander--Reiten $\E$-triangle
$\delta$.  Fix $X\in\ind\C$.  The projective boundary argument in the proof
of Theorem~\ref{thm:main}, using enough projectives and the local finiteness
of $\Pcal$, shows that \[
 \{P\in\ind\Pcal\mid\C(P,X)\ne0\}
\]
is finite.  Choose $Q_X$ as in {(D$_-$)} and an $\E$-triangle
\[
 \Omega Q_X\longrightarrow P_X\longrightarrow Q_X\dashrightarrow,
~~\mbox{where}~P_X\in\Pcal.
\]
Since its relation belongs to $\Rel(\C)=\AR(\C)$, it is a finite integral
combination of Auslander--Reiten relations. By Proposition~\ref{prop:Serre-EFS}
and Lemma~\ref{lem:stable-defect}, we imply that $ \{U\in\ind\C\mid\underline\C(U,Q_X)\ne0\}$
is finite.  Condition {(D$_-$)} now gives
$|\Supp\C(-,X)|<\infty$.

Since $D$ is a duality on finite-length $R$-modules, it reflects the zero
module.  Serre duality therefore yields, for indecomposable $X$ and $U$,
\[
 \C(X,U)\ne0
 \quad\Longleftrightarrow\quad
 \C(U,SX)\ne0.
\]
Consequently $ \Supp\C(X,-)=\Supp\C(-,SX),$
as sets of isomorphism classes.  Since $S$ is an autoequivalence, $SX$ is indecomposable, and the right
side is finite by the first part of the proof.  Thus $\C$ is locally finite.
\end{proof}

The Serre functor in Theorem~\ref{thm:Serre-one-sided} is defined on $\C$
itself.  It should not be confused with Auslander--Reiten--Serre duality for
extriangulated categories.  The latter is formulated on the stable and
costable categories, see \cite{IyamaNakaokaPalu2024}.  In particular, an
identity
$$
 D\E(X,Y)\simeq\C(Y,TX)
$$
with $T$ an autoequivalence of $\C$ forces every projective and every
injective object to be zero.  Indeed, if $P$ is projective, then
$0=D\E(P,TP)\simeq\End(TP)$.  Hence $TP=0$, and therefore $P=0$.
The injective case is dual.

\vspace{5mm}

\hspace{-5.5mm}{\bf Statement on AI}\hspace{2mm} ChatGPT was used only for language editing and improving the clarity of the manuscript. All mathematical results, arguments, and proofs were developed and verified by the authors.
\vspace{3mm}

\hspace{-5.5mm}\textbf{Data Availability}\hspace{2mm} Data sharing not applicable to this article as no datasets were generated or analysed during
the current study.
\vspace{3mm}

\hspace{-5.5mm}\textbf{Conflict of Interests}\hspace{2mm} The authors declare that they have no conflicts of interest to this work.

\hspace{-6.5mm}
\textbf{Panyue Zhou}\\
School of Mathematics and Statistics, Changsha University of Science and Technology, 410114 Changsha, Hunan, P. R. China\\
Email: panyuezhou@163.com\\

\hspace{-5mm}\textbf{Bin Zhu}\\
Department of Mathematical Sciences, Tsinghua University, 100084 Beijing, P. R. China\\
Email: zhu-b@mail.tsinghua.edu.cn.

\end{document}